\documentclass{article}
\usepackage[utf8]{inputenc}
\usepackage[a4paper, total={15.92cm, 24.37cm}, top=1in, left=1in]{geometry}
\usepackage[T1]{fontenc}
\usepackage{lmodern}
\usepackage{dsfont} \newcommand{\ds}{\mathds} \newcommand{\cl}{\mathcal}
\usepackage{graphicx}
\usepackage{amsmath}
\usepackage{amsthm}
\newtheorem{thm}{Theorem}
\newtheorem{prp}{Proposition}[section]
\newtheorem{dfn}[prp]{Definition}
\newtheorem{lmm}[prp]{Lemma}
\newtheorem{cvn}[prp]{Convention}
\newtheorem{rmk}[prp]{Remark}

\newcommand{\R}{\ensuremath{\mathds{R}}}
\renewcommand{\P}{\ensuremath{\mathds{P}}}

\newcommand{\var}{\ensuremath{\text{\rm Var}}}			

\newcommand{\lp}{\ensuremath{\left(}}
\newcommand{\rp}{\ensuremath{\right)}}
\newcommand{\lb}{\ensuremath{\left\{}}
\newcommand{\rb}{\ensuremath{\right\}}}
\newcommand{\ls}{\ensuremath{\left[}}
\newcommand{\rs}{\ensuremath{\right]}}

\newcommand{\LL}{\ensuremath{\mathcal{L}}}
\newcommand{\XX}{\ensuremath{\mathcal{X}}}

\newcommand{\DD}{\ensuremath{\mathcal{D}}}

\renewcommand{\SS}{\ensuremath{\mathcal{S}}}

\newcommand{\RR}{\ensuremath{\mathcal{R}}}

\title{%
    On soft capacities, quasi-stationary distributions\\
    and the pathwise approach to metastability
}
\author{%
    A. Bianchi\thanks{%
        Università degli Studi di Padova,
        Dipartimento di Matematica,
        Via Trieste, 63, 35121 Padova,
        Italy,
        e-mail: bianchi@math.unipd.it
    }
    \and
    A. Gaudillière\thanks{%
        Aix Marseille Univ, CNRS, Centrale Marseille, I2M,
        Marseille,
        France,
        e-mail: alexandre.gaudilliere@math.cnrs.fr
    }
    \and
    P. Milanesi\thanks{%
        Aix Marseille Univ, CNRS, Centrale Marseille, I2M,
        Marseille,
        France,
        e-mail: paolo.milanesi@univ-amu.fr
    }
}
\date{}

\begin{document}
\maketitle
\begin{abstract}
    Motivated by the study of the metastable stochastic Ising model
    at subcritical temperature
    and in the limit of a vanishing magnetic field,
    we extend the notion of $(\kappa, \lambda)$-capacities
    between sets, as well as the associated notion of soft-measures,
    to the case of {\it overlapping} sets.
    We recover their essential properties,
    sometimes in a stronger form
    or in a simpler way,
    relying on weaker hypotheses.
    These properties allow to write
    the main quantities associated
    with reversible metastable dynamics,
    e.g.\ asymptotic transition and relaxation times,
    in terms of objects
    that are associated  with two-sided variational principles.
    We also clarify the connection with the classical ``pathwise approach''
    by referring to temporal means on the appropriate time scale.

    \par\smallskip\noindent
    {\it Acknowledgements:} A. G. and P. M. thank Maria Eulalia Vares,
    the Universidade federal do Rio de Janeiro and the Università di Padova
    for the kind hospitality
    which gave us the opportunity to lay the foundations of this work.
\end{abstract}

\section{Model and main results}\label{mercurio}

In the present paper we consider
a general Markovian model
for a metastable dynamics
and we show how,
under mild hypotheses,
``soft measures and capacities''
associated with a covering
of the configuration space
allow for a description of the metastable state
provide sharp estimates of the relaxation time,
fast relaxation to local equilibria,
and enable to establish the asymptotic exponential law
of the transition time to equilibrium.
The ultimate goal of this paper is to provide
a mathematical framework
to prove a convergence in law
to an exponential distribution
of the transition time to equilibrium
for the metastable kinetic Ising model
at any subcritical temperature
and in the limit of a vanishing magnetic field.
In the companion paper~\cite{GMV},
we establish such a convergence
by working out the model-dependent part of the proof.
In addition, we present an explicit
connection with the classical ``pathwise approach to metastability'',
and we provide a comparison with the techniques
appeared in the recent literature of abstract metastable dynamics.

\begin{figure}
\caption{Snapshots of a kinetic Ising model
at temperature $T = 1.5$
and under magnetic field $h = 0.14$
in a $256 \times 256$ box
with periodic boundary conditions.
Minus spins are yellow and plus spins are red.
We took pictures 
at times 471, 7482 and 13403 for the first line,
at times 14674, 15194, 15432, 15892 and 16558 for the second line,
and times 17328, 23645 and 40048 for the last line.
}\label{cacio_e_pepe}
$$
\includegraphics[clip, trim = {1mm 1mm 1mm 1mm}, width = 1.8 cm]{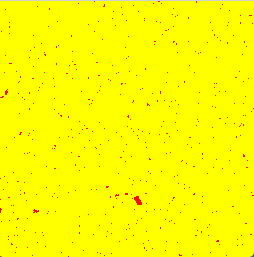}
\qquad
\includegraphics[clip, trim = {1mm 1mm 1mm 1mm}, width = 1.8 cm]{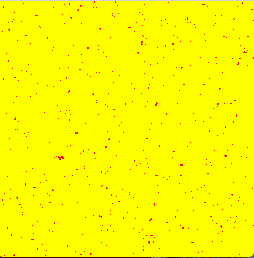}
\qquad
\includegraphics[clip, trim = {1mm 1mm 1mm 1mm}, width = 1.8 cm]{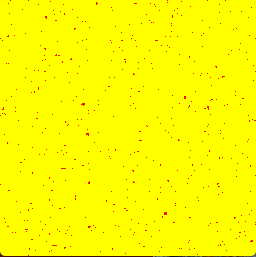}
$$
$$
\includegraphics[clip, trim = {1mm 1mm 1mm 1mm}, width = 1.8 cm]{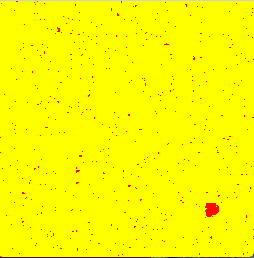}
\qquad
\includegraphics[clip, trim = {1mm 1mm 1mm 1mm}, width = 1.8 cm]{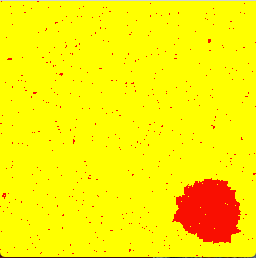}
\qquad
\includegraphics[clip, trim = {1mm 1mm 1mm 1mm}, width = 1.8 cm]{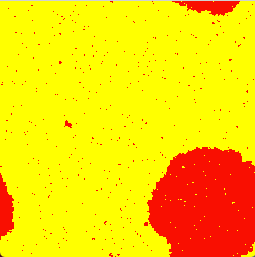}
\qquad
\includegraphics[clip, trim = {1mm 1mm 1mm 1mm}, width = 1.8 cm]{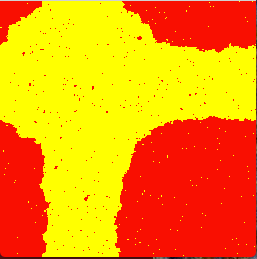}
\qquad
\includegraphics[clip, trim = {1mm 1mm 1mm 1mm}, width = 1.8 cm]{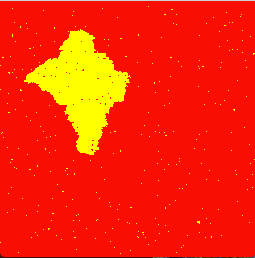}
$$
$$
\includegraphics[clip, trim = {1mm 1mm 1mm 1mm}, width = 1.8 cm]{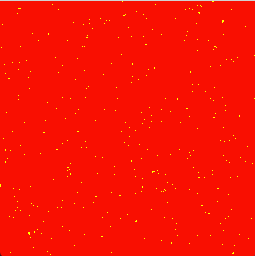}
\qquad
\includegraphics[clip, trim = {1mm 1mm 1mm 1mm}, width = 1.8 cm]{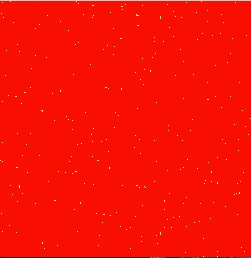}
\qquad
\includegraphics[clip, trim = {1mm 1mm 1mm 1mm}, width = 1.8 cm]{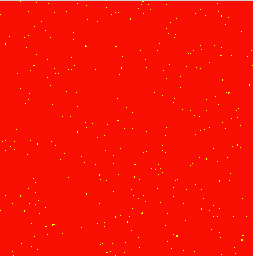}
$$
\end{figure}

Before stating the main results,
and to better illustrate the ideas leading to the present research,
we start by discussing the qualitative behavior of a concrete Markovian model.
We consider a continuous time Glauber dynamics
in two-dimensional finite boxes
of area diverging as $1 / h^2$,
where $h \ll 1$ is the magnetic field, as studied in~\cite{GMV}.
By Glauber dynamics we mean a single spin flip dynamics
that is reversible with respect to the Gibbs measure
of the Ising model, thus including Metropolis and
heat-bath dynamics (see~\cite{GMV} for precise definitions).
Figure~\ref{cacio_e_pepe} shows a sample
of such a Metropolis dynamics
started from  the ``all minus configuration'',
with periodic boundary conditions, and with
an Hamiltonian that is twice\footnote{%
    This is because~\cite{GMV} sticks to the conventions
    of its main reference~\cite{SS},
    where the authors introduced
    a relatively unusual factor $1 / 2$ in the Hamiltonian.
    Since Figure~\ref{cacio_e_pepe}
    illustrates a dynamics ran at inverse temperature
    $\beta = 2 / 3$ without such a convention,
    it would correspond to a trajectory sampled at inverse
    temperature $\beta' = 4 / 3$,
    i.e., to a twice smaller temperature,
    with the convention of~\cite{SS} and~\cite{GMV}.
}
the Hamiltonian of~\cite{GMV},
up to boundary conditions.

\goodbreak
Figure~\ref{cacio_e_pepe} is the kind of picture
to bear in mind when reading the present paper.
Here are some of its most relevant features.
\begin{itemize}
\item The last three pictures can be seen as different samples
    from the equilibrium Gibbs measure.
    Since $h > 0$, it is concentrated on configurations
    where plus spins dominate.
    One has to wait for some ``transition to equilibrium''
    to start observing non-anomalous configurations
    with respect to this equilibrium state.
\item The first three pictures can be seen as different samples
    from a single distribution concentrated on configurations
    where minus spins prevail.
    Notice that they represent typical configurations
    of a ``local equilibrium'', or metastable state,
    that is very different  from the equilibrium stable state.
\item The transition to equilibrium
    is triggered by the nucleation of a ``supercritical droplet''
    (in the first picture, see the tiny plus-phase droplet
    that was not large enough to trigger it).
    Notice that the time needed to invade the whole box
    is of a smaller order than the time needed to appear
    ---about 2000 and 15000 time units, respectively---
    as a fluctuation of the metastable state.
\item This suggests that the time needed by the system
    to relax to some local equilibrium (metastable or stable state)
    is short with respect to the global relaxation time.
    This is what we call ``fast relaxation to local equilibria''.
\item This also suggests that the transition to equilibrium
    can be described as ``both late and abrupt'',
    and it is thus compatible
    with an asymptotic exponential law of the transition time to equilibrium.
\end{itemize}

Topping the list of questions posed when trying to build
a rigorous mathematical counterpart of such observations:
\begin{itemize}
\item[i.] How can we make sense of such a notion of
    local or metastable equilibrium?
    How should it be described?
\item[ii.] How can we define such a random transition time to equilibrium?
    How should its law be described?
\end{itemize}
The pathwise approach to metastability introduced in~\cite{CGOV}
and fully described in~\cite{OV} suggests using temporal means
and stopping times to answer these questions.
If
$$
    X: t \geq 0 \longmapsto X(t) \in \cl X
$$
is the Markov process, the metastable character of which we want to investigate,
and $f: \cl X \rightarrow \ds R$ is an observable,
let us denote the associated temporal mean from time $t \geq 0$
and on time scale $\theta > 0$ by
\begin{equation}\label{antoine}
    A_\theta(t, f) = {1 \over \theta} \int_t^{t + \theta} f(X(s))\,ds.
\end{equation}
Let $\mu$ be the equilibrium distribution
of the process and, for any probability measure $\nu$
on $\cl X$ and any observable $f$,
let us denote by $\nu(f)$ the mean value of $f$ with respect to $\nu$.
Within the pathwise approach, answering the previous questions
amounts to identifying a probability measure $\tilde \mu$,
a stopping time $T$ and a deterministic time scale $\theta$
such that, for the system started in $\tilde \mu$:
\begin{itemize}
\item $\theta$ is small with respect to $\ds E_{\tilde \mu}[T]$,
    the mean value of $T$;
\item $A_\theta(t, f)$ is close to $\tilde \mu(f)$
    for all $t < T - \theta$ and all observable $f$
    with very large probability;
\item $A_\theta(t, f)$ is close to $\mu(f)$ for all observables $f$
    and all $t \geq T$ of the same order as $\ds E_{\tilde \mu}[T]$,
    or possibly larger, with very large probability.
\end{itemize}
The metastable state and the transition time to equilibrium
can then be defined by $\tilde \mu$ and $T$, respectively.
The pathwise approach also focuses on establishing the convergence
in law to an exponential distribution with parameter one
of the rescaled transition time $T / \ds E_{\tilde \mu}[T]$.
This is to account for a ``both late and abrupt''
transition to equilibrium,
by opposition with the behaviour of a slowly relaxing system
for which this relaxation would occur by following
a gentle and constant drift\footnote{%
    In the case of the kinetic Ising model,
    one can instead think of two opposite strong drifts:
    one towards the metastable state and one towards the stable state.
    The nucleation of the supercritical droplet occurs by fluctuation
    {\it against\/} the first drift
    and the system follows the second drift afterwards.
    This is coherent with the fast relaxation to local equilibria.
}.
One cannot distinguish between these two behaviours
by looking only at the time evolution of the mean value
$$
    t \geq 0 \longmapsto \ds E_{\tilde \mu}\bigl[f(X(t))\bigr],
    \qquad
    f: \cl X \longrightarrow \ds R,
$$
that are associated with the natural Markovian semi-group;
hence the terminology ``{\it pathwise} approach''.

As far as the construction of the probability measure $\tilde \mu$
and the stopping time $T$ is concerned,
the authors of~\cite{CGOV} suggested,
by following~\cite{LP},
to use a {\it restricted ensemble\/}
$$
    \mu_{\cl R} = \mu\bigl(\cdot \bigm| \cl R\bigr)
$$
and the associated exit time $T_{\cl X \setminus \cl R}$,
for some subset $\cl R$ of the configuration space $\cl X$
that can be thought of as a ``basin of attraction''
of the metastable state.
In the case of the kinetic Ising model,
Schonmann and Sholsman made a detailed study
in the beautiful paper~\cite{SS}
of such restricted ensemble and exit time
by defining $\cl R$ as the set of configurations
for which the + spins can be enclosed into small enough
(subcritical) contours
(see~\cite{SS} or the companion paper~\cite{GMV}
for more details).
Their results allow to prove that the temporal means
$A_\theta(t, f)$, for a suitable $\theta$,
are close to $\mu_{\cl R}(f)$ or $\mu(f)$,
but only for times $t$
that are very small or very large
with respect to $\ds E_{\mu_{\cl R}}[T_{\cl X \setminus \cl R}]$,
which is strongly related with the lack of control
of the law of the transition time to equilibrium.

In~\cite{BG} we proposed some modifications
of the traditional construction of the metastable state $\tilde \mu$
and the transition time to equilibrium $T$.
We worked with quasi-stationary measures
to connect the two issues
of the description of the metastable state
and of the law of the transition time.
This led us first to a mathematical formalization of
the phenomenon of ``local thermalization''
or ``fast relaxation to local equilibria'',
then to concieve of the stable state 
as a more stable local equilibrium only.
As a consequence, we suggested
to use a stopping time at which one is ensured\footnote{%
    Note that for the kinetic Ising model,
    despite the strong drift towards equilibrium
    at the appearance of the first supercritical droplet,
    there is still some time to wait before the relaxation
    to be achieved:
    the pictures on the second line
    of Figure~\ref{cacio_e_pepe} show
    atypical configurations with respect to the equilibrium measure.
    The exit from $\cl R$ does not coincide with global thermalization.
} to be close
to this more stable equilibrium,
instead of the exit time $T_{\cl X \setminus \cl R}$.
This eventually led us to the notion of soft measure,
which is a quasi-stationary measure (of a trace process)
associated with such a stopping time rather
than a plain exit time:
while quasi-stationary measures
are associated with absorption
as soon as a Markov process
hits some set $\cl S$,
which acts as an {\it hard\/} barrier,
we consider absorption, or killing,
at some rate $\lambda$ in $\cl S$
---our Markov process can then penetrate $\cl S$,
hence the terminology ``{\it soft\/} measure''.
We showed that, with a suitable set of hypotheses
in which $\cl S$ can be envisaged
as a basin of attraction of the stable state,
and for a suitable choice of the killing parameter $\lambda$,
the associated soft measure $\mu^*_{\cl R, \lambda_{\cl S}}$
and killing time $T_{\lambda_{\cl S}}$
(see sections~\ref{julie} and~\ref{novara} for precise definitions)
were natural candidates to build the metastable state $\tilde \mu$
and the transition time $T$.
This came with an associated notion of ``soft capacity''
and a solution,
as detailed in~\cite{BG} and shortly reviewed
in Section~\ref{clo} of the present paper,
to difficulties encountered in proving
asymptotic exponential laws for transition times
and that had also emerged
in the ``potential-theoretic approach'' to metastability.

Unfortunately, we could not prove the convergence in law
of a rescaled transition time to an exponential random variable
in the aforementioned case of the kinetic Ising model
within the framework of~\cite{BG}.
The problem came from the fact that, in~\cite{BG},
we assumed that the sets $\cl R$ and $\cl S$ form
a partition of the configuration space.
But as shown in~\cite{GMV}, the tools
to control local relaxation times,
that is the crucial point
of the model-dependent part
of the study,
involve some mesoscopic and macroscopic analysis
that is incompatible with such a neat separation
at the microscopic level.
It turned out that removing the disjointness hypothesis
of $\cl R$ and $\cl S$ from the framework of~\cite{BG},
with the aim of improving
the global proposed strategy
and  eventually prove an asymptotic exponential law
in~\cite{GMV}, required new ideas.
This was the first motivation for
writing the present paper.
In doing so, we introduced some simplifications and
strengthened some estimates with respect to~\cite{BG},
and we also provided an explicit connection
with the classical pathwise approach.

At last, we have to make a final remark
on the asymptotic exponential law for the kinetic Ising model.
We prove such an asymptotic result in the companion paper~\cite{GMV},
where we consider the kinetic Ising model
in a finite box of area diverging as $1 / h^2$,
in the limit of a vanishing magnetic field $h > 0$.
Even though the proof of such a result was missing,
this was to be expected and nothing to be surprised at.
But in the case of infinite volume dynamics,
as the one considered in~\cite{SS}, the situation is more delicate,
and it is not so clear whether an asymptotic exponential law
for the transition time to equilibrium
is the correct conjecture for the both late and abrupt transition.
Indeed, as shown by Shonmann and Shlosman,
the transition mechanism of infinite volume dynamics
does not only involve the nucleation of supercritical droplets,
which should occur after an asymptotic exponential time,
but also the coalescence of such growing droplets
on a time scale that is logarithmically equivalent
to the nucleation time scale.
Thus, to decide whether the asymptotic exponential law
still holds when the times involved in both mechanisms are sequentially added,
a precise control on the prefactors of their mean values is needed.
This underlines the relevance of the sharp estimates
of spectral gap and mean transition time
(see Theorem~\ref{albero} in Section~\ref{julie})
in terms of soft capacities inherited from
the classical ones that were introduced by~\cite{BEGK02},
\cite{BEGK04} and~\cite{BEGK05} in metastability studies.

The paper is organized as follows.
In Section~\ref{julie} we present our main result,
in Section~\ref{clo} we discuss the relations with the more recent literature
on abstract metastable dynamics,
in Section~\ref{sole} we collect all the material
that can be directly imported from~\cite{BG}
and we simplify some of its proofs,
and in sections~\ref{luna}, \ref{pranzo} and~\ref{telefono}
we give the proofs of our main theorems.

\subsection{Local and global thermalizations}\label{julie}

Let $X$ be an irreducible continuous time Markov process
on a finite configuration space
$$
    \cl X = \cl R \cup \cl S
$$
with $\cl R$ and $\cl S$ possibly overlapping
and with generator ${\cal L}$ defined by
\begin{equation}\label{venere}
    ({\cal L}f)(x) = \sum_{y \neq x} w(x, y) [f(y) - f(x)],
    \qquad x \in {\cal X},
    \qquad f: {\cal X} \rightarrow \mathds{R}.
\end{equation}
We assume $X$ to be reversible
with respect to a probability measure $\mu$
and we denote by $\cl D$ the Dirichlet form
associating with any observable $f : \cl X \rightarrow \ds R$
$$
    \cl D(f)
    = {1 \over 2} \sum_{x, y \in \cl X} \mu(x)w(x, y)\bigl[
        f(x) - f(y)
    \bigr]^2.
$$
We call $1 / \gamma$ the relaxation time of the dynamics,
so that $\gamma$ is the spectral gap
\begin{equation}\label{reggio_emilia}
    \gamma
    = \min_{
        \scriptstyle f \in \ell^2(\mu),
    \atop
        \scriptstyle {\rm Var}_{\mu}(f) \neq 0
    } {{\cal D}(f) \over {\rm Var}_{\mu}(f)}
    \,.
\end{equation}
With any $\cl Y \subset \cl X$
we associate the restricted dynamics $X_{\cl Y}$ in $\cl Y$
with generator defined by
\begin{equation}\label{marte}
    ({\cal L}_{\cl Y}f)(x)
    = \sum_{y \in \cl Y \setminus \{x\}} w(x, y) [f(y) - f(x)],
    \qquad x \in {\cal Y},
    \qquad f: {\cal Y} \rightarrow \mathds{R},
\end{equation}
which is reversible with respect to the restricted ensemble
$$
    \mu_{\cl Y} = \mu\bigl(\cdot \bigm| \cl Y\bigr).
$$
If $X_{\cl Y}$ is also an irreducible process
we denote by $\gamma_{\cl Y}$ its spectral gap.
The associated exit rate (see Section~\ref{sole}
for an explanation of the terminology) is
$$
    \phi^*_{\cl Y}
    = \min_{
        \scriptstyle f \neq 0,
    \atop
        \scriptstyle f|_{{\cal X}\backslash{\cal Y}} = 0
    } {{\cal D}(f) \over E_\mu[f^2]}
    \,.
$$
Let $\sigma_\lambda$ be an exponential random variable
of parameter $\lambda$, independent of $X$, and
$\ell_{\cl S}(t)$ denote the local time of $X$ in $\cl S$, that is
$$
    \ell_{\cl S}(t) = \int_0^t \ds 1_{\{X(s) \in \cl S\}} \, ds\,.
$$
We then define $T_{\lambda_{\cl S}}$ as the first time $t \geq 0$
such that $\ell_{\cl S}(t)$ reaches $\sigma_\lambda$:
$$
    T_{\lambda_{\cl S}} = \min\bigl\{
        t \geq 0 : \ell_{\cl S}(t) \geq \sigma_\lambda
    \bigr\}\,.
$$
In other words, $T_{\lambda_{\cl S}}$ 
is the absorption time of the process $X$
killed at rate $\lambda_{\cl S}$
with 
\begin{equation}\label{pierre}
    \lambda_{\cl S}(x) = \lambda \ds 1_{\{x \in \cl S\}},
    \qquad x \in \cl X.
\end{equation}
We define analogously $T_{\kappa_{\cl R}}$,
which is associated with the killing rate $\kappa_{\cl R}$,
from an exponential random variable $\sigma_\kappa$
of parameter $\kappa$, independent of $X$ and $\sigma_\lambda$,
and by using $\ell_{\cl R}$, the local time in $\cl R$.
With
$$
    V_\kappa^\lambda(x) = \mathds{P}_x\left(
        T_{\kappa_{\cal R}} < T_{\lambda_{\cal S}}
    \right),
    \qquad x \in {\cal X},
$$
the associated soft capacity, or $(\kappa, \lambda)$-capacity is
$$
    C_\kappa^\lambda({\cal R}, {\cal S})
    = {\cal D}\bigl(V_\kappa^\lambda\bigr)
    + \kappa \mu({\cal R}) E_{\mu_{\cal R}}\left[
        \left(
            V_\kappa^\lambda - 1
        \right)^2
    \right]
    + \lambda \mu({\cal S}) E_{\mu_{\cal S}}\left[
        \left(
            V_\kappa^\lambda - 0
        \right)^2
    \right].
$$

We assume to be in some asymptotic regime,
meaning that the dynamics depends on some parameters,
as temperature, volume or magnetic field
going to some finite or infinite limit,
and we will write $f \ll g$ and $f \sim g$
whenever $0 \leq f \leq g$ and the ratio $f / g$
goes to zero and one, respectively,
in the considered asymptotic regime.
To state the first main result of this paper
we introduce
\medskip\par\noindent
{\bf Hypotheses $(H)$:} \label{falchetto} {\it
    We assume ${\cal R}$ and ${\cal S}$
    to be such that
    \begin{itemize}
    \item[1.] the associated restricted processes
    $X_{\cal R}$, $X_{{\cal R} \backslash {\cal S}}$,
    $X_{\cal S}$ and $X_{{\cal S} \backslash {\cal R}}$
    are all irreducible;
\item[2.]  $\phi_{{\cal R} \backslash {\cal S}}^* \ll \gamma_{\cal R}$
    and $\phi_{{\cal S} \backslash {\cal R}}^* \ll \gamma_{\cal S}$;
\item[3.]
    $\mu({\cal S}) \geq \mu({\cal R)}$ and
    $\phi_{{\cal R} \backslash {\cal S}}^* \ll \gamma_{\cal S}$.\end{itemize}}

\medskip\par\noindent
Since exit rates $\phi^*_{\cl Y}$ are associated
with some mean exit time from $\cl Y$
(see Section~\ref{novara}),
such hypotheses amount in particular to say
that ``local relaxation times'' $1 / \gamma_{\cl Y}$
are small with respect to these mean exit times.
To check $(H)$ we need upper bounds on exit rates
and lower bounds on spectral gaps.
While the former are usually quite easy to get,
since exit times are given by an infimum
in some variational principle,
getting the latter is the core of the difficulty.
This is the crux of the companion paper~\cite{GMV}.
However, if the hypothesis set $(H)$ is in force,
even from a rough control of local relaxation times
we can get sharp estimates of the global one:
\begin{thm}\label{albero}
    Assuming hypotheses $(H)$ and choosing $\kappa$ and $\lambda$
    such that
    $$
        \phi^*_{{\cal R} \backslash {\cal S}}
        \ll \kappa
        \ll \gamma_{\cal R}
        \qquad {\it and} \qquad
        \phi^*_{{\cal R} \backslash {\cal S}}
        \vee \phi^*_{{\cal S} \backslash {\cal R}}
        \ll \lambda
        \ll \gamma_{\cal S},
    $$
    it holds
    $$
        {1 \over \gamma} \sim {
            \mu({\cal R}) \mu({\cal S})
            \over
            C_{\kappa}^{\lambda}({\cal R}, {\cal S})
        }
        \qquad {\it and} \qquad
        \ds E_{\mu_{\cl R}}\bigr[T_{\lambda_{\cl S}}\bigl]
        \sim {
            \mu({\cal R})
        \over
            C_{\kappa}^{\lambda}({\cal R}, {\cal S})
        }\,.
    $$
\end{thm}

\goodbreak

\noindent {\bf Comments:}
\begin{itemize}
\item[i.] This gives a way to compute the precise asymptotic value
    of the relaxation time $1 / \gamma$. As soon as $\kappa$ and $\lambda$
    are chosen in the right window, the spectral gap $\gamma$ scales,
    up to the multiplicative term $\mu({\cal R})\mu({\cal S})$,
    like the soft capacity $C_\kappa^\lambda({\cal R}, {\cal S})$
    that satisfies a two-sided variational principle.
    As we will see in Section~\ref{sole},
    any test function and flow
    will give upper and lower bound
    for $C_\kappa^\lambda({\cal R}, {\cal S})$.
    These are then translated in upper and lower bounds for $\gamma$.
\item[ii.] We already mentioned that the main
    difficulty in proving the validity of the hypothesis set $(H)$
    usually lies in bounding $\gamma_{\cal R}$ and $\gamma_{\cal S}$
    from below. It is often the case that,
    once such lower bounds are established,
    there is a natural way to build a suitable flow
    to get useful lower bounds of the soft capacity:
    see \cite{BG} and \cite{GMV}.
\item[iii.] We will see with the next theorem
    that $T_{\lambda_{\cl S}}$ can be identified,
    under slightly stronger hypotheses,
    as a transition time to equilibrium if $\mu(\cl R) \ll \mu(\cl S)$,
    or to a more stable local equilibrium
    in the general case $\mu(\cl R) \leq \mu(\cl S)$.
    Theorem~\ref{albero} says in particular that,
    when starting from the restricted ensemble $\mu_{\cl R}$,
    its expected value and the relaxation time
    have the same asymptotic behaviour
    up to the multiplicative factor $1 / \mu(\cl S)$, which lies between~1
    and~2 in the latter case and goes to~1 in the former one.
    To prove this result we will focus on the local time in $\cl R$
    up to transition
    $$
        T^{\cl R}_{\lambda_{\cl S}}
        = \ell_{\cl R}\bigl(T_{\lambda_{\cl S}}\bigr)
    $$
    and its expected value
    $$
        {1 \over \phi^*_{\cl R, \lambda_{\cl S}}}
        = \ds E_{\mu^*_{\cl R, \lambda_{\cl S}}}\bigl[
            T^{\cl R}_{\lambda_{\cl S}}
        \bigr]
    $$
    for the dynamics started from the soft measure
    $\mu^*_{\cl R, \lambda_{\cl S}}$ defined in Section~\ref{novara}.
    For $\lambda$ ranging from~0 to $+\infty$ the soft measure family
    provides a continuous interpolation
    from the restricted ensemble $\mu_{\cl R}$
    to the quasi-stationary distribution $\mu^*_{\cl R \setminus \cl S}$
    associated with the process immediately killed when hitting $\cl S$.
    We will see in Section~\ref{sole} that hypotheses $(H)$
    implies that our soft measures are all close together,
    each of them providing a metastable state model.
    And we will prove that
    $$
        E_{\mu_{\cl R}}\bigr[T_{\lambda_{\cl S}}\bigl]
        \sim {1 \over \phi^*_{\cl R, {\lambda_{\cl S}}}}
    $$
    for any $\lambda$ in the allowed window of Theorem~\ref{albero}.
\item[iv.] The proof of Theorem~\ref{albero} is deferred to Section~\ref{luna}.
    In this section we will give explicit upper and lower bounds
    for $\phi^*_{{\cal R}, \lambda_{\cal S}}$ and $\gamma$
    that will bound explicitly the multiplicative corrections
    with respect to the asymptotic values.
\item[v.] About these upper and lower bounds,
    the latter are the most relevant.
    Both the exit rate $\phi^*_{{\cal R}, \lambda_{\cl S}}$
    and the spectral gap $\gamma$ are indeed {\it infima}
    given by some variational principles,
    so that upper bound are simpler to get.
    In particular, the lower bound given in Section~\ref{luna}
    on $\gamma$ is a Poincaré inequality
    that can be generalized to the case
    of more than one metastable state
    (cf. Proposition~\ref{giggiola}).
    The upper bounds indicate however that our lower bounds
    catch the right orders.
\item[vi.] The windows given in Theorem~\ref{albero}
    for $\kappa$ and $\lambda$ can in principle be enlarged
    for the conclusion to hold.
    For example, as far as the lower bound on $\kappa$
    is concerned, the explicit upper and lower bound of Section~\ref{luna}
    will scale like the given asymptotic values
    for $\phi^*_{{\cal R}, \lambda_{\cal S}} / \kappa \ll 1$.
    We will see in Section~\ref{sole} that
    $\kappa \gg \phi^*_{{\cal R} \backslash {\cal S}}$
    is only a sufficient condition for having
    $\phi^*_{{\cal R}, \lambda_{\cal S}} / \kappa \ll 1$.
    Our hypotheses $(H)$ can be relaxed in the same way
    for the conclusion to hold.
\item[vii.] When going from the case ${\cal R} \cap {\cal S} = \emptyset$
    of \cite{BG} to the case ${\cal R} \cap {\cal S} \neq \emptyset$
    of the present paper,
    some of the arguments of \cite{BG} are immediately generalized,
    while others require more work.
    The new arguments of Section~\ref{luna} not only allow for weaker hypotheses,
    they also lead to a better upper bound for $\phi^*_{\cl R, \lambda_{\cl S}}$,
    already in the special case $\cl S = \cl X \setminus \cl R$.
\end{itemize}

\goodbreak

To connect hypotheses $(H)$ with the pathwise approach framework,
we set $T_0 = 0$, $T_1 = T_{\lambda_{\cl S}}$,
we denote by $T_2$
the absorption time associated with the killing rate
$\kappa_{\cl R}$ after $T_1$, that is
$$
   T_2 = T_1 + T_{\kappa_{\cal R}} \circ \Theta_{T_1}
$$
---where $\Theta_t$ is the usual shift operator, such that, for all $s \geq 0$,
$(\Theta_t (X))(s) = X(t + s)$--- and we set for all $i \geq 0$,
\begin{align*}
    T_{2i + 1} &= T_{2i} + T_{\lambda_{\cal S}} \circ \Theta_{T_{2i}},\\
    T_{2i + 2} &= T_{2i + 1} + T_{\kappa_{\cal R}} \circ \Theta_{T_{2i + 1}},
\end{align*}
using at each step independent copies of the exponential variables
$\sigma_\kappa$ and $\sigma_\lambda$.
Since we will work with mixing times
rather than relaxation times associated with restricted processes,
we introduce
\begin{equation}\label{uovo}
    \chi_{\cal R} = \max_{x \in {\cal R}} {1 \over \mu_{\cal R}(x)}\,,
    \qquad
    \chi_{\cal S} = \max_{x \in {\cal S}} {1 \over \mu_{\cal S}(x)}\,,
\end{equation}
and a strengthened hypothesis set.
\medskip\par\noindent
{\bf Hypotheses $(H')$:} {\it
    Hypotheses $(H)$ hold together with
    $(\ln \chi_{\cal S}) / \gamma_{\cal S}
    \ll 1 / \phi^*_{{\cal R} \backslash {\cal S}}$
    and
    $(\ln \chi_{\cal R}) / \gamma_{\cal R}
    \ll 1 / \phi^*_{{\cal R} \backslash {\cal S}}$.
}
\medskip\par\noindent
Recalling Equation~\eqref{antoine} it holds:
\begin{thm}\label{corteo}
    Assuming hypotheses (H') and choosing $\kappa$ and $\lambda$ such that
    $$
        \phi^*_{{\cal R} \backslash {\cal S}}
        \ll \kappa
        \ll \gamma_{\cal R} / \ln \chi_{\cal R}
        \qquad
        \hbox{and}
        \qquad
        \phi^*_{{\cal R} \backslash {\cal S}}
        \ll \lambda
        \ll \gamma_{\cal S} / \ln \chi_{\cal S},
    $$
    there is a positive $\theta \ll 1 / \phi^*_{{\cal R}, \lambda_{\cal S}}$
    such that, for all $\delta > 0$, for all $n > 0$
    and for all observable $f : {\cal X} \rightarrow \mathds{R}$,
    $$
        \mathds{P}_{\mu_{\cal R}}\left(
            \begin{array}{l}
                \forall i < n,\;
                T_{2i + 1} > T_{2i} + \theta,\;
                T_{2i + 2} > T_{2i + 1} + \theta,\\[4pt]
                \forall t \in [T_{2i},  T_{2i + 1} - \theta[,\;
                \forall t' \in [T_{2i + 1}, T_{2i + 2} - \theta[,\\[4pt]
                |A_\theta(t, f) - \mu_{\cal R}(f)| \leq \delta,\;
                |A_\theta(t', f) - \mu_{\cal S}(f)| \leq \delta
            \end{array}
        \right) \sim 1.
    $$
\end{thm}

\goodbreak
\noindent {\bf Comments:}
\begin{itemize}
\item[i.] We will give the proof of this theorem in Section~\ref{telefono}.
    We will make explicit
    a possible choice for $\theta$
    and we will give explicit estimates
    of such a probability.
\item[ii.] Let us consider the empirical distribution process
    defined by
    $$
        A_\theta(t) = {1 \over \theta} \int_t^{t + \theta} \delta_{X(s)}\,ds,
        \qquad t \geq 0\,.
    $$
    If one is interested in convergence of $A_\theta$,
    one can consider (as in~\cite{OV}, Chapter~4)
    a slight modification of the process $A_\theta(t)$.
    With the above notation, let us consider the modification
    $$\tilde A_\theta(t) = \left\{
        \begin{array}{ll}
            \tilde A_\theta(T_k - \theta) &
            \mbox{ if } t\in [T_k - \theta, T_k[ \mbox{ for some }k,\\
              A_\theta(t) &
             \mbox{ otherwise.}
         \end{array}\right.\,.
    $$
    Then Theorem \ref{corteo},
    together with a convergence in law of the rescaled $T_k$,
    provides the key tool to prove convergence
    in Skorokhod topology of $\tilde A_\theta$.
\item[iii.] In the case $\mu({\cal R}) \ll \mu({\cal S})$,
    we will have (recall Theorem~\ref{albero} and its associated third comment)
    $$
        {1 \over \gamma}
        \sim \ds E_{\mu_{\cl R}}\bigl[
            T_{\lambda_{\cl S}}
        \bigr]
        \sim {1 \over \phi^*_{{\cal R}, \lambda_{\cal S}}}
        \ll {1 \over \phi^*_{{\cal S}, \kappa_{\cal R}}}
        \sim \ds E_{\mu_{\cl S}}\bigr[
            T_{\kappa_{\cl R}}
        \bigl]
    $$
    with $\phi^*_{\cl S, \kappa_{\cl R}}$ symmetrically defined
    with respect to $\phi^*_{\cl R, \lambda_{\cl S}}$.
    The rescaled process on time unit $1 / \phi^*_{{\cal R}, \lambda_{\cal S}}$
    reaches equilibrium
    after time $\phi^*_{{\cal R}, \lambda_{\cal S}} T_1$ of order 1,
    this equilibrium can be described by $\mu_{\cal S}$,
    and there is no oscillation on this time scale between the two states
    $\mu_{\cal R}$ and $\mu_{\cal S}$.
    Our system will go back to $\mu_{\cal R}$
    only after a much longer time of order
    $1 / \phi^*_{{\cal S}, \kappa_{\cal R}}$.
\end{itemize}

As far as we are concerned with the asymptotic exponential law
of $T_{\lambda_{\cl S}}$ and $T_{\kappa_{\cl R}}$,
when starting from $\mu_{\cl R}$ and $\mu_{\cl S}$, respectively,
we can adapt and even simplify the arguments of~\cite{BG}.
This is a simple consequence of the following facts,
that are the contents of Lemma~\ref{saturno},
Lemma~\ref{neptuno} and Proposition~\ref{monza}.
(See also Proposition~\ref{gazza} for a quantitative statement.)
\begin{itemize}
\item[i.] $T^{\cl R}_{\lambda_{\cl S}}
     = \ell_{\cal R}\bigl(T_{\lambda_{\cal S}}\bigr)$
    is an exact exponential random variable,
    of rate $\phi^*_{\cl R, \lambda_{\cl S}}$,
    when starting from $\mu^*_{\cl R, \lambda_{\cl S}}$.
    This is the reason to work with such quasi-stationary distributions.
\item[ii.] Since $\ell_{\cl S \setminus \cl R}\bigl(T_{\lambda_{\cl S}}\bigr)$
    is dominated by the exponential random variable $\sigma_\lambda$,
    taking $\lambda$ in the right window and passing to the expected value,
    it holds that
    $$
        {1 \over \lambda}
        \ll {1 \over \phi^*_{\cl R \setminus \cl S}}
        \leq {1 \over \phi^*_{\cl R, \lambda_{\cl S}}}.
    $$
    Hence, since  $T_{\lambda_{\cl S}}
    = \ell_{\cl R}\bigl(T_{\lambda_{\cl S}}\bigr)
    + \ell_{\cl S \setminus \cl R}\bigl(T_{\lambda_{\cl S}}\bigr)$,
    and together with point i., we deduce the asymptotic law
    of $T_{\lambda_{\cl S}}$
    when starting from $\mu^*_{\cl R, \lambda_{\cl S}}$.
\item[iii.] The restricted ensemble $\mu_{\cl R}$
    and the soft measure $\mu^*_{\cl R, \lambda_{\cl S}}$
    are close in total variation distance.
\end{itemize}

Our last main result is a ``local thermalization property''
which implies the existence of an asymptotic exponential law
for stopping times $T_{\lambda_{\cal S}}$ or $T_{\kappa_{\cal R}}$,
independently of the choice of the initial distribution.
As in \cite{BG} we can use
$T_{\lambda_{\cal S}}$ and $T_{\kappa_{\cal R}}$
to build a stopping time $T^*$
of order $1 / (\kappa \wedge \lambda)$ at most
and such that the law of $X(T^*)$ is close
to $\mu^*_{{\cal R}, \lambda_{\cal S}}$
or $\mu^*_{{\cal S}, \kappa_{\cal R}}$
by conditioning on two complementary events.
Since we are here mainly interested
in convergence in law of the transitions times,
closeness in total variation distance
will be enough.
We can then build $T^*$ in a simpler way,
and with weaker hypotheses,
than in \cite{BG}.

We recall that the previously mentioned
and symmetric exit rates $\phi^*_{\cl R, \lambda_{\cl S}}$
and $\phi^*_{\cl S, \kappa_{\cl R}}$
are precisely defined in Section~\ref{novara}.
Their asymptotic behaviour is precisely controlled
in terms of soft capacities by Theorem~\ref{pisa}
in Section~\ref{brune}.
\begin{thm}\label{gheppio}
    Assuming hypotheses (H') and choosing
    $\kappa$ and $\lambda$ such that
    $$
        \phi^*_{{\cal R} \backslash {\cal S}}
        \ll \kappa
        \ll \gamma_{\cal R} / \ln\chi_{\cal R}
        \qquad\hbox{and}\qquad
        \phi^*_{{\cal R} \backslash {\cal S}}
        \ll \lambda
        \ll \gamma_{\cal S} / \ln\chi_{\cal S},
    $$
    there are three stopping times $T^*_{\cal R}$,
    $T^*_{\cal S}$ and $T^* = T^*_{\cal R} \wedge T^*_{\cal S}$
    such that, whatever the starting distribution of~$X$:
    \begin{itemize}
    \item[i.] the expected value of $T^*$
        is smaller than or equal to $2 / (\kappa \wedge \lambda)$;
    \item[ii.]
        with $\pi^*_{\cl R}$ and $\pi^*_{\cl S}$
        the laws of $X(T^*)$
        conditioned to $\{T^* = T^*_{\cal R}\}$
        and $\{T^* = T^*_{\cal S}\}$, respectively,
        it holds
        $$
            d_{\rm TV}\bigl(
                \mu_{\cal R}, \pi^*_{\cl R}
            \bigr) \leq {3 \kappa \over \gamma_{\cal R}}\left(
                1 + \left[
                    \ln {\gamma_{\cal R} \sqrt{\chi_{\cal R}} \over 2 \kappa}
                \right]_+
            \right) \ll 1,
        $$
        and
        $$
            d_{\rm TV}\bigl(
                \mu_{\cal S}, \pi^*_{\cl S}
            \bigr) \leq {3 \lambda \over \gamma_{\cal S}}\left(
                1 + \left[
                    \ln {\gamma_{\cal S} \sqrt{\chi_{\cal S}} \over 2 \lambda}
                \right]_+
            \right) \ll 1;
        $$
    \item[iii.]
        conditionally to $\{T^* = T^*_{\cal R}\}$ or $\{T^* = T^*_{\cal S}\}$,
        the random variable
        $\phi^*_{{\cal R}, \lambda_{\cal S}}
        T_{\lambda_{\cal S}} \circ \Theta_{T^*}$
        or
        $\phi^*_{{\cal S}, \kappa_{\cal R}}
        T_{\kappa_{\cal R}} \circ \Theta_{T^*}$, respectively,
        converge in law to an exponential variable or parameter~1;
    \item[iv.] it holds
        $$
            \mathds{P}_\nu\left(
                T^* = T^*_{\cal R} = T_{\kappa_{\cal R}} < T_{\lambda_{\cal S}}
            \right)
            \geq 1 - 3 \mathds{P}_\nu\left(
                T_{\lambda_{\cal S}} < T_{\kappa_{\cal R}}
            \right)
        $$
        and
        $$
            \mathds{P}_\nu\left(
                T^* = T^*_{\cal S} = T_{\kappa_{\cal S}} < T_{\lambda_{\cal R}}
            \right)
            \geq 1 - 3 \mathds{P}_\nu\left(
                T_{\lambda_{\cal R}} < T_{\kappa_{\cal S}}
            \right),
        $$
        so that these probabilities converge to~1
        when
        $\mathds{P}_\nu(T_{\lambda_{\cal S}} < T_{\kappa_{\cal R}}) \ll 1$
        and
        $\mathds{P}_\nu(T_{\kappa_{\cal R}} < T_{\lambda_{\cal S}}) \ll 1$,
        respectively.
    \end{itemize}
\end{thm}

\goodbreak
\noindent {\bf Comments:}
\begin{itemize}
\item[i.] We will give the proof of this theorem in Section~\ref{pranzo},
    together with the explicit construction of $T^*$
    and explicit estimates on all error terms.
\item[ii.] This theorem implies an asymptotic exponential law
    for the rescaled time
    $\phi^*_{{\cal R}, \lambda_{\cal S}} T_{\lambda_{\cal S}}$
    when $X$ starts from any measure $\nu$
    such that $\mathds{P}_\nu(T_{\lambda_{\cal S}} < T_{\kappa_{\cal R}}) \ll 1$.
    We will give a direct proof of this fact in Section~\ref{sole}
    (see Proposition~\ref{tortora}).
    One of the reasons why we needed stronger hypotheses in \cite{BG}
    was that we were also interested there
    in events like $\mathds{P}_\nu(T_{\lambda_{\cal S}} > t)$
    for times $t$ that were large with respect to $1 / \gamma_{\cal R}$
    but of a smaller order than $1 / \phi^*_{{\cal R}, \lambda_{\cal S}}$.
\item[iii.] Theorem~\ref{gheppio} implies an analogue of Theorem~\ref{corteo}
    where $T_0$ is replaced by $T^*$,
    $T_1$ by $T^* + T_{\lambda_{\cal S}} \circ \Theta_{T^*}$
    and $\mathds{P}_{\mu_{\cal R}}$
    by $\mathds{P}_x\bigl(\, \cdot \bigm| T^* = T^*_{\cal R}\bigr)$,
    for any $x$ in ${\cal X}$.
\end{itemize}

\subsection{The entropy issue}\label{clo}

The problem of describing metastable states
goes back to~\cite{Max}.
But we will limit our discussion of the related literature
to the last five decades by starting from~\cite{LP}.
Since such a review,
albeit already biased,
was developed and detailed in~\cite{BG},
we tried to make it short in the present paper.

After the introduction in~\cite{LP} of restricted ensembles
for describing metastable states,
the pathwise approach introduced in~\cite{CGOV}
was mainly built on~\cite{FW} and its large deviations techniques.
It entangled the first two main questions
---that of describing the metastable state
and that of characterizing the transition time to equilibrium---
with a third one, that of giving a detailed picture
of the typical paths followed along the transition\footnote{%
    In the case of our kinetic Ising model
    this amounts to describing
    the shape of subcritical, critical and subcritical droplets.
    The techniques presented in~\cite{GMV} allow
    for proving that the critical droplet is Wulff-shaped
    and this question is still largely open
    for subcritical and supercritical droplets.
}.
One of the main difficulties faced when dealing
with each particular model was to give {\it upper\/} bounds
for the expected value of the transition time to equilibrium.
In the case of ``very low temperature'' dynamics
or ``rare transition limits'', this question
is related with that of giving a detailed description
of the ``(virtual) energy landscape''.
Both reversible and non-reversible dynamics were considered
within the pathwise approach,
we refer to~\cite{OV} for a full picture
and to~\cite{CNS} for what we believe to be
the most detailed and general account,
along Freidlin and Wentzell's line,
on this main difficulty
(disentangled from the third question)
in the case of very low temperature systems.

Potential-theoretic tools introduced by~\cite{BEGK02},
\cite{BEGK04} and~\cite{BEGK05} in metastability studies
brought in particular two major improvements:
\begin{itemize}
\item  they related the estimation of mean transition times
    with capacity estimates for which variational principles were available,
    without reference to any energy landscape;
\item the upper and lower bounds obtained from these variational principles,
    led in many cases to sharp estimates of mean transition times.
\end{itemize}
Basic hypotheses expressed in terms of mean hitting times
and capacities allowed for asymptotic exponential laws
for transition times between local equilibria.
We refer to~\cite{BdH} for a full account on these techniques,
initially restricted to the reversible case.
These potential-theoretic tools were incorporated
in~\cite{Bl10} and~\cite{BL11},
where the authors proposed
to deal with metastable dynamics by proving,
using martingale equations,
a convergence in law to a Markovian dynamics in Skorohod topology
of the label evolution associated with a trace process.
From~\cite{GL} and~\cite{Lan} to~\cite{LMS},
this was extended to non-reversible dynamics.

When we wrote the first version of~\cite{BG},
the vast majority of metastability studies
were restricted either to very low temperature dynamics,
or to Markovian models that could be naturally mapped
to dynamics with very low entropy metastable states,
such as distributions that are
concentrated around a single point of the configuration space.
There were two kinds of difficulties
when dealing with metastable states
for which entropy really matters:
local relaxation was difficult to prove
and there was no clear road map to follow.
This was the case for the kinetic Ising model
at {\it fixed} subcritical temperature,
that is the entropy issue.
As far as the general theory was concerned,
when looking for example at potential-theoretic tools
that linked the estimation of mean transition times to capacities
(for which variational formulas were available)
one came across two difficulties among others:
the formula used harmonic measures that were not
naturally linked with metastable states,
and it required estimates
for a mean electrostatic potential
that did not satisfy any variational formula.
It seemed also that one could get asymptotic exponential laws
only if one was able to solve these difficulties.
By introducing soft measures and capacities,
we reduced in~\cite{BG} these difficulties of the general theory
to a single model-dependant issue: getting (rough) upper bounds
for the relaxation times of restricted dynamics.
Theorem~\ref{albero} of the present paper indeed provides
an asymptotic formula for the mean transition time
started in the metastable state,
and not from a harmonic measure\footnote{%
    This is however an asymptotic formula
    and there is no equality.
}, that does not require
any estimate on the mean potential $\mu(V_\kappa^\lambda)$.
We transposed the difficulties of the general theory
to that of checking in each particular model
the hypotheses $(H)$ or $(H')$, namely to getting
upper bounds for relaxation times of restricted dynamics.
Such solutions were immediately adapted to the framework
of the ``martingale approach to metastability''~\cite{BL15}.
Armendáriz,  Grosskinsky and Loulakis
could then study the metastable condensing zero-range process
in the thermodynamic limit~\cite{AGL}.

We already mentioned that we could not directly apply
the results of~\cite{BG}
to deal with the kinetic Ising model:
\cite{BG} requires $\cl R$ and $\cl S$ to form a partition
of the configuration space,
but the mesoscopic and macroscopic
analysis we make in~\cite{GMV} to control the relaxation times
of the restricted dynamics does not allow
for such a neat separation at the microscopic level.
In the martingale approach framework,
local equilibria are however associated
with usually far apart subsets $\cl R_i$ of the configuration space.
This space is written as a disjoint union
$$
    \cl X = \cl N \cup \bigcup_{i \in I} \cl R_i,
$$
where $\cl N$ is a subset of $\cl X$
in which the system is expected to spend a negligible time
on the time scale that allows for the transition to equilibrium.
With such a framework,
we might still be able to bound the relaxation times
on the restricted dynamics in $\cl R_i$
(say taking $\cl R_0 \subset \cl R$ and $\cl R_1 \subset \cl S$,
with $\cl R$ and $\cl S$ defined in~\cite{GMV})
and we could avoid the issue of configurations
which in~\cite{GMV} belong to both $\cl R$ and $\cl S$:
such critical configurations would have been in $\cl N$.
But this would have only displaced some difficulties.
As previously mentioned,
one of the main issues is that of getting upper bounds
for the mean transition time.
This is turned by Theorem~\ref{albero},
and in the martingale approach,
into a problem of lower bounds for a (soft) capacity.
Even though associated with a variational principle,
this is often one of the most delicate point.
As mentioned in the second comment of Theorem~\ref{albero},
such lower bounds are however obtained simply, in~\cite{GMV},
from the control of the local relaxation times.
In getting such lower bounds for the soft capacities,
we do use the fact that $\cl R$ and $\cl S$
are ``neighbouring'' sets (even overlapping in this case).
Using overlapping $\cl R$ and $\cl S$, instead of far apart $\cl R_i$
together with a negligible set $\cl N$ can be a matter of taste.
But there is an issue when estimating the relevant capacities.

It is worth noting at this point
that in going from restricted ensemble
to soft measures, we start changing our representation for the metastable set,
by considering overlapping $\cl R$ and $\cl S$ we stop considering only
local equilibria with disjoint support,
and our soft capacities,
that could be considered as our main objects, now only depend
on some killing rates like $\lambda_{\cal S}$ and not $\lambda$ and $\cl S$
separately.
We could actually consider capacities
associated with full support killing rates,
and much more general local equilibria.
Similar ideas have been considered in~\cite{ACGM17}
with concrete applications~\cite{ACGM20} in a quite different field.

As far as the exponential law in itself is concerned
\cite{FMNS}, \cite{BG} and~\cite{FMNSS}
contain similar detailed and quantitative results,
obtained with rather similar techniques.
The first paper is limited to low entropy metastable states
and uses recurrence to an associated reference configuration,
the second paper explains how this recurrence property
can be replaced by a fast relaxation hypothesis
to a quasi-stationary measure
and such is the approach of the third paper.
By contrast we only focus in the present paper
on convergence in law results,
which are still quantitative but less precise
and associated with simpler constructions.
But as in~\cite{BG} we also study, on the one hand,
the convergence to a stable equilibrium
and not only the exit time of a domain,
which is the main object of interest
in~\cite{FMNS} and~\cite{FMNSS},
and we are interested, on the other hand,
in getting sharp bounds
for mean transition times
by relying on capacity estimates:
this is where our reversibility hypothesis matters,
while~\cite{FMNS} and~\cite{FMNSS}
do not assume reversibility.

We close this discussion on the recent literature
of abstract metastable dynamics by listing,
in order of importance,
the improvement of the present paper
with respect to~\cite{BG}:
\begin{itemize}
\item By relaxing the previous hypothesis $\cl R \cap \cl S = \emptyset$,
    we provided a framework to prove, in the companion paper~\cite{GMV},
    the asymptotic exponential law of the two-dimensional kinetic Ising model,
    at any subcritical temperature and in the limit of vanishing magnetic field;
\item we made explicit the connection with the classical pathwise
    approach by referring to temporal means and identifying explicitly
    an associated time scale $\theta$;
\item we simplified the construction of the short time scale stopping time $T^*$
    at which the system is distributed
    ``as $\mu_{\cl R}$ or $\mu_{\cl S}$'' independently of  the starting measure;
\item we simplified the proof of the asymptotic exponential law;
\item we obtained a sharper quantitative control of the mean transition time
    to equilibrium.
\end{itemize}

\section{Basic definitions and properties}
\label{sole}

\subsection{Soft measures and killed process}\label{novara}
Let $X$ be an irreducible Markov process,
with generator ${\cal L}$ as in Equation~\eqref{venere}.
We assume $X$ to be reversible with respect to a probability measure $\mu$.
For ${\cal R}$ and ${\cal S}$
such that ${\cal X} = {\cal R} \cup {\cal S}$, and
for any given finite $\lambda \geq 0$
we build in this section
the soft measure $\mu_{{\cal R}, \lambda_{\cal S}}^*$
from the trace $X^*_{{\cal R}, \lambda_{\cal S}}$ on ${\cal R}$
of the process $X$ killed at rate $\lambda_{\cal S}(x)$
in each $x$ in ${\cal X}$.
To facilitate the reading,
we use the upper index $^*$ to denote objects associated with killed processes.

Formally, let $T_{\cal R}$ be the hitting time of ${\cal R}$,
and $T_{\lambda_{\cal S}}$ be the absorption time of $X$
killed at rate $\lambda_{\cl S}$ defined in Equation~\eqref{pierre}.
Equivalently\footnote{%
    This is an alternative construction
    that does not use the exponential variable
    $\sigma_\lambda$ of Section~\ref{julie}
    and allows for less homogeneous killing rates.
}, we can associate any $x\in X$ with an independent Poisson process
of rate $\lambda_{\cal S}(x)$, and
define $T_{\lambda_{\cal S}}$ as the first arrival time $t$
of the Poisson process associated with $x = X(t)$.
The trace process $X^*_{{\cal R}, \lambda_{\cal S}}$
is then the Markov process on ${\cal R}$
that, on one hand, is killed in $x$ in ${\cal R}$
at rate
\begin{equation}\label{mantova}
    e^*_{{\cal R}, \lambda_{\cal S}}(x)
    = \lambda_{\cal S}(x) + \sum_{z \neq x} w(x, z) \mathds{P}_z\left(
        T_{\cal R} > T_{\lambda_{\cal S}}
    \right)
    = \lambda_{\cal S}(x) + \sum_{z \not\in {\cal R}} w(x, z) \mathds{P}_z\left(
        T_{\cal R} > T_{\lambda_{\cal S}}
    \right),
\end{equation}
where the second term accounts for the rate
with which $X$ escapes from ${\cal R}$
and is killed before returning in it,
and that, on the other hand, jumps from $x$ to a distinct $y \in {\cal R}$
at rate
\begin{equation}\label{cremona}
    w_{{\cal R}, \lambda_{\cal S}}(x, y)
    = \sum_{z \neq x} w(x, z) \mathds{P}_z\left(
        X(T_{\cal R}) = y, T_{\cal R} < T_{\lambda_{\cal S}}
    \right)
    = w(x, y)  + \sum_{z \not\in {\cal R}} w(x, z) \mathds{P}_z\left(
        X(T_{\cal R}) = y, T_{\cal R} < T_{\lambda_{\cal S}}
    \right).
\end{equation}
It is associated with the sub-Markovian generator
${\cal L}^*_{{\cal R}, \lambda_{\cal S}}$
that acts on functions $f: {\cal R} \rightarrow \mathds{R}$ according to
\begin{equation}\label{generatore}
    \left(
        {\cal L}^*_{{\cal R}, \lambda_{\cal S}} f
    \right)(x)
    = -e^*_{{\cal R}, \lambda_{\cal S}}(x) f(x)
    + \sum_{y \in {\cal R} \backslash \{x\}}
        w_{\cal R, \lambda_{\cal S}}(x, y) \bigl(f(y) - f(x)\bigr),
    \qquad x \in {\cal R}.
\end{equation}

Notice that the irreducibility hypothesis (cf. $(H)$) on $X_R$,
having generator given by Equation~\eqref{marte},
implies the irreducibility of $X^*_{{\cal R}, \lambda_{\cal S}}$
associated with Equation~\eqref{generatore}.
Then, by Perron-Frobenius Theorem,
${\cal L}^*_{{\cal R}, \lambda_{\cal S}}$
has only negative eigenvalues
and the smallest eigenvalue $\phi^*_{{\cal R}, \lambda_{\cal S}} > 0$
of $-{\cal L}^*_{{\cal R}, \lambda_{\cal S}}$ is non-degenerate
and associated with a left eigenvector $\mu^*_{{\cal R}, \lambda_{\cal S}}$
which is a probability on $\cal R$.
We call it the {\it soft measure}
associated with the killing rate $\lambda_{\cal S}$
and it is a special kind of quasi-stationary measure.
We denote by $\ell_{\cal R}$ the local time in ${\cal R}$:
$$
    \ell_{\cal R}(t) = \int_0^t \mathds{1}_{\{X(u) \in {\cal R}\}} \, du,
    \qquad t \geq 0.
$$
It is a non-decreasing continuous function of the time $t$
and we call $\ell_{\cal R}^{-1}$ its right-continuous inverse:
$$
    \ell^{-1}_{\cal R}(s)=\inf \left\{t\geq 0 : \ell_{\cal R}(t) > s\right\},
    \qquad s \geq 0.
$$
We then have
$X^*_{{\cal R}, \lambda_{\cal S}}(s) = X \circ \ell_{\cal R}^{-1}(s)$
for all
$s < T^{\cal R}_{\lambda_{\cal S}}
= \ell_{\cl R}\bigl(T_{\lambda_{\cl S}}\bigr)$
and the following lemma.
\begin{lmm}\label{saturno}
    For all $s \geq 0$ and all $x$ and $y$ in ${\cal R}$ it holds
    $$
        \mathds{P}_{\mu^*_{\RR,\lambda_\SS}}\lp
            X\circ\ell_{\RR}^{-1}(s) = x
            \bigm| T_{\lambda_\SS}^{\cal R}>s
        \rp = \mu_{\RR,\lambda_{\cal S}}^*(x),
    $$
    $$
        \lim_{s \to \infty} \P_x\lp
            X\circ\ell_{\RR}^{-1}(s) = y
            \bigm| T_{\lambda_\SS}^{\cal R} > s
        \rp = \mu_{\RR,\lambda_{\cal S}}^*(y),
    $$
    and
    $$
        \P_{\mu_{\RR,\lambda_{\cal S}}^*}\lp
            T_{\lambda_\SS}^{\RR} > s
        \rp = e^{-s\phi_{\RR,\lambda_{\cal S}}^*}.
    $$
    The exit rate is also given by
    $$
        \phi_{\RR,\lambda_{\cal S}}^*
        = \mu_{\RR,\lambda_{\cal S}}^*\lp
            e^*_{\RR,\lambda_{\cal S}}
        \rp.
    $$
\end{lmm}
These formulas are consequences of the fact that,
for all $s \geq 0$,
every function $f$ and every probability $\nu$ on ${\cal R}$,
$$
    \mathds{E}_{\nu}\left[
        f\bigl(X^*_{{\cal R}, \lambda_{\cal S}}(s)\bigr)
        \mathds{1}_{\left\{
            T_{\lambda_{\cal S}}^{\cal R} > s
        \right\}}
    \right] = \nu \left(
        e^{s {\cal L}^*_{{\cal R}, \lambda_{\cal S}}} f
    \right).
$$
We refer to Lemma~2.12 in \cite{BG}
for more details.

As in \cite{BG} (Lemma 2.13) and assuming hypotheses $(H)$,
when $\lambda$ goes from $0$ to $+ \infty$
we get a continuous interpolation between
the restricted ensemble $\mu_{\cal R}$
and the quasi-stationary measure $\mu^*_{{\cal R} \backslash {\cal S}}$.
Indeed, when $\lambda = 0$ our soft measure $\mu^*_{{\cal R}, \lambda_{\cal S}}$
is nothing but the invariant measure of the trace process on ${\cal R}$,
which coincides with $\mu(\cdot | {\cal R})$
as a consequence of our reversibility hypothesis.
And, if $\lambda = +\infty$ we are led to define $X^*_{{\cal R}, \lambda_{\cal S}}$
as $X^*_{{\cal R} \backslash {\cal S}}$,
obtained from $X$ by killing it instantaneously in ${\cal S}$.
This process is associated
with the sub-Markovian generator ${\cal L}^*_{{\cal R} \backslash {\cal S}}$
defined by
$$
    \left(
        {\cal L}^*_{{\cal R} \backslash {\cal S}} f
    \right)(x)
    = -e^*_{{\cal R} \backslash {\cal S}}(x) f(x)
    + \sum_{
        \scriptstyle y \in {\cal R} \backslash {\cal S},
    \atop
        \scriptstyle y \neq x
    } w(x, y) (f(y) - f(x)),
    \qquad f: {\cal R} \backslash {\cal S} \rightarrow \mathds{R},
    \qquad x \in {\cal R},
$$
with, for all $x$ in ${\cal R} \backslash {\cal S}$,
\begin{equation}\label{urano}
    e^*_{{\cal R} \backslash {\cal S}}(x) = \sum_{z \in {\cal S}} w(x, z).
\end{equation}
Since $X^*_{{\cal R} \backslash {\cal S}}$ is assumed
to be irreducible, the quasi-stationary measure $\mu^*_{{\cal R} \backslash {\cal S}}$
is unambiguously defined.
Contrarily to $X^*_{{\cal R} \backslash {\cal S}}$,
which cannot penetrate ${\cal S}$,
$\mu^*_{{\cal R}, \lambda_{\cal S}}$ and $X^*_{{\cal R}, \lambda_{\cal S}}$
are associated with a process which can visit the killing zone ${\cal S}$.
As mentioned earlier,
this is the explanation for the terminology {\it soft} measure.

We denote by $\|\cdot\|$ and $\langle \cdot, \cdot \rangle$
the standard norm and scalar product in $\ell^2(\mu)$.
For any ${\cal A} \subset {\cal X}$
we denote by $\mu_{\cal A}$ the conditional probability measure
$\mu_{\cal A} = \mu(\cdot | {\cal A})$,
and by $\|\cdot\|_{\cal A}$ and $\langle\cdot, \cdot\rangle_{\cal A}$
the standard norm and scalar product in
$\ell^2(\mu_{\cal A})$.
We also denote by $e^*_{\cal A}$, as in Equation~\eqref{urano},
the escape rate
\begin{equation*}
    e^*_{\cal A}(x) = \sum_{z \not\in {\cal A}} w(x, z)
\end{equation*}
for $x$ in ${\cal A}$.
\begin{prp}\label{milano}
    When $\lambda$ grows to $+\infty$,
    \begin{equation} \label{giove}
        \phi^*_{{\cal R}, \lambda_{\cal S}}
        = \min_{
                \scriptstyle f: {\cal R} \rightarrow \mathds{R}
        \atop
                \scriptstyle f \neq 0
        }{
            \langle f, -{\cal L}^*_{{\cal R}, \lambda_{\cal S}} f\rangle_{\cal R}
        \over
            \|f\|^2_{\cal R}
        }
        = \min_{f > 0} {
            \langle f, -{\cal L}^*_{{\cal R}, \lambda_{\cal S}} f\rangle_{\cal R}
        \over
            \|f\|^2_{\cal R}
        }
    \end{equation}
    continuously grows to $\phi^*_{{\cal R} \backslash {\cal S}}$.
    In particular, for all $\lambda \geq 0$, it holds
    \begin{equation}\label{plutone}
        \phi^*_{{\cal R}, \lambda_{\cal S}}
        \leq \phi^*_{{\cal R} \backslash {\cal S}}
        = \min_{
            \scriptstyle f: {\cal R} \backslash {\cal S} \rightarrow \mathds{R}
        \atop
            \scriptstyle f \neq 0
        }{
            \langle
                f, -{\cal L}^*_{{\cal R} \backslash {\cal S}} f
            \rangle_{{\cal R} \backslash {\cal S}}
        \over
            \|f\|_{{\cal R} \backslash {\cal S}}^2
        }
        \,.
    \end{equation}
\end{prp}
\noindent {\it Proof:\/}
The first equality in Equation~\eqref{giove}
comes from the fact that ${\cal L}^*_{{\cal R}, \lambda_{\cal S}}$
is a self-adjoint operator.
The last equality in Equation~\eqref{giove}
is a consequence of Perron-Frobenius Theorem.
Continuity follows from Lemma~\ref{saturno}
and continuity properties of $\mu^*_{{\cal R}, \lambda_{\cal S}}$
and  $e^*_{{\cal R}, \lambda_{\cal S}}$.
Finally, for any positive function $f$ on ${\cal R}$
it holds
\begin{align*}
    \langle f, -{\cal L}^*_{{\cal R}, \lambda_{\cal S}} f \rangle_{\cal R}
    &= - \sum_{x \in {\cal R}} \sum_{y \in {\cal R} \backslash \{x\}} \mu_{\cal R}(x) f(x) w_{{\cal R}, \lambda_{\cal S}}(x, y) f(y) \\
    &\quad + \sum_{x \in {\cal R}} \mu_{\cal R}(x) \left(
        e^*_{{\cal R}, \lambda_{\cal S}}(x) + \sum_{y \in {\cal R} \backslash \{x\}}  w_{{\cal R}, \lambda_{\cal S}}(x, y)
    \right) f^2(x).
\end{align*}
The first term, minus sign included,
which accounts for the cross product,
is clearly increasing in $\lambda$,
since $w_{{\cal R}, \lambda_{\cal S}}(x, y)$
decreases with $\lambda$ and $f$ is positive.
As far as the second term is concerned, using the definitions (\ref{mantova}) and \ref{cremona}, we compute
\begin{align*}
    e^*_{{\cal R}, \lambda_{\cal S}}(x)
    + \sum_{y \in {\cal R} \backslash \{x\}}  w_{{\cal R}, \lambda_{\cal S}}(x, y)
    &= \lambda_{\cal S}(x) + \sum_{z \not\in {\cal R}} w(x, z) \mathds{P}_z\left(
        T_{\cal R} > T_{\lambda_{\cal S}}
    \right) \\
    &\quad + \sum_{y \in {\cal R} \backslash \{x\}}  w(x, y)
    + \sum_{y \in {\cal R} \backslash \{x\}} \sum_{z \not\in {\cal R}} w(x, z) \mathds{P}_z\left(
        X(T_{\cal R}) = y, T_{\cal R} < T_{\lambda_{\cal S}}
    \right)\\
    &= \lambda_{\cal S}(x) + \sum_{z \not\in {\cal R}} w(x, z) \mathds{P}_z\left(
        T_{\cal R} > T_{\lambda_{\cal S}}
    \right)
    + \sum_{y \in {\cal R} \backslash \{x\}}  w(x, y)\\
    &\quad + \sum_{z \not\in {\cal R}} w(x, z) \Bigl(
        1 - \mathds{P}_z\left(
            T_{\cal R} > T_{\lambda_{\cal S}}
        \right)
        - \mathds{P}_z \left(
            X(T_{\cal R}) = x, T_{\cal R} < T_{\lambda_{\cal S}}
        \right)
    \Bigr) \\
    &= \lambda_{\cal S}(x) +  w(x)
     - \sum_{z \not\in {\cal R}} w(x, z) \mathds{P}_z\left(
            X(T_{\cal R}) = x, T_{\cal R} < T_{\lambda_{\cal S}}
        \right)\,,.
\end{align*}
where in the last line we set $w(x)\,=\, \sum_{z\neq x} w(x,z)$.
The whole coefficient is then increasing in $\lambda$
and the monotonicity of $\phi^*_{{\cal R}, \lambda_{\cal S}}$
follows.
\qed

Using simple test functions in such variational principles
we get the following upper bounds on $\phi^*_{{\cal R}, \lambda_{\cal S}}$.
\begin{lmm}\label{neptuno}
    It holds, for all $\lambda \geq 0$,
    \begin{equation}\label{neptuno_1}
        \phi^*_{{\cal R}, \lambda_{\cal S}}
        \leq \phi^*_{{\cal R} \backslash {\cal S}}
        \leq \mu_{{\cal R} \backslash {\cal S}}\bigl(e^*_{{\cal R} \backslash {\cal S}}\bigr),
    \end{equation}
    \begin{equation}\label{neptuno_2}
        \phi^*_{{\cal R}, \lambda_{\cal S}}
        \leq \lambda \mu_{\cal R}({\cal R} \cap {\cal S}) + \mu_{\cal R}\bigl(e^*_{\cal R}\bigr),
    \end{equation}
    \begin{equation}\label{neptuno_3}
        \mathds{E}_{\mu_{\cal R}}\left[
            T^{\cal R}_{\lambda_{\cal S}}
        \right]
        \leq {1 \over \phi^*_{{\cal R}, \lambda_{\cal S}}}
        \,.
    \end{equation}
\end{lmm}
\noindent{\it Proof:\/}
Formula~\eqref{neptuno_1} is obtained from Formula~\eqref{plutone}
with the constant function equal to $1$ on ${\cal R} \backslash {\cal S}$
used as a test function.
We use the constant function $\mathds{1}_{\cal R}$ equal to $1$ on ${\cal R}$
in Equation~\eqref{giove} and we bound by $1$
the last probability appearing in Equation~\eqref{mantova}
to get Inequality~\eqref{neptuno_2}.
Finally, to get Inequality~\eqref{neptuno_3},
we observe that $1 / \phi^*_{{\cal R}, \lambda_{\cal S}}$
is the largest eigenvalue of the Green kernel
${\cal G}^R_{\lambda_{\cal S}} = \bigl(- {\cal L}_{{\cal R}, \lambda_{\cal S}}\bigr)^{-1}$,
so that,
writing $\ell_y$ for the local time in $y$,
\par
\vbox{$$
    {1 \over \phi^*_{{\cal R}, \lambda_{\cal S}}}
    = \max_{\|f\|_{\cal R} = 1} \langle
        f, {\cal G}^{\cal R}_{\lambda_{\cal S}}f
    \rangle_{\cal R}
    \geq \langle
        \mathds{1}_{\cal R}, {\cal G}^{\cal R}_{\lambda_{\cal S}}\mathds{1}_{\cal R}
    \rangle_{\cal R}
    = \sum_{x, y \in {\cal R}} \mu_{\cal R}(x) \mathds{E}_x\left[
        \ell_y\bigl(T_{\lambda_{\cal S}}\bigr)
    \right]
    = \mathds{E}_{\mu_{\cal R}}\left[
        \ell_{\cal R}\bigl(T_{\lambda_{\cal S}}\bigr)
    \right].
$$
\qed}
\par\noindent{\bf Remark:}
While $\phi^*_{{\cal R}, \lambda_{\cal S}}$
and $\phi^*_{{\cal R} \backslash {\cal S}}$
are the means of escape rates according to
quasi-stationary measures that can be difficult to compute,
the upper bounds given in \eqref{neptuno_1} and \eqref{neptuno_2}
refer only to the reversible measure $\mu$.

\subsection{Restricted dynamics}\label{poiana}
For any finite $\lambda \geq 0$,
the {\it restricted dynamics} $X_{{\cal R}, \lambda_{\cal S}}$
is obtained from $X^*_{{\cal R}, \lambda_{\cal S}}$ or $X$
by removing all killing moves or killing excursions outside ${\cal R}$.
It is associated with the generator ${\cal L}_{{\cal R}, \lambda_{\cal S}}$
that acts on functions $f: {\cal R} \rightarrow \mathds{R}$ according to
$$
    \left(
        {\cal L}_{{\cal R}, \lambda_{\cal S}} f
    \right)(x)
    = \sum_{y \in {\cal R} \backslash \{x\}} w_{\cal R, \lambda_{\cal S}}(x, y) (f(y) - f(x)),
    \qquad x \in {\cal R},
$$
with $w_{{\cal R}, \lambda_{\cal S}}(x, y)$
defined by Equation~\eqref{cremona}
for distinct $x$ and $y$ in ${\cal R}$.
In particular this process in reversible with respect
to the restricted measure $\mu_{\cal R}$
and it inherits irreducibility from the restricted process $X_{\cal R}$.
We call $\gamma_{{\cal R}, \lambda_{\cal S}}$ its spectral gap
\begin{equation}\label{lodi}
    \gamma_{\RR,\lambda_\SS}
    =\min_{\substack{f:\RR\to\mathds{R}, \\ \var_{\mu_{\RR}}\lp f\rp\neq 0}}\frac{\DD_{\RR,\lambda_\SS}(f)}{\var_{\mu_{\RR}}\lp f\rp}
    \,,
\end{equation}
with ${\cal D}_{{\cal R}, \lambda_{\cal S}}$ the {\it Dirichlet form}
defined by
$$
    \DD_{\RR,\lambda_\SS}(f)
    = \langle
        f, -{\cal L}_{{\cal R}, \lambda_{\cal S}}f
    \rangle_{\cal R}
    =\frac{1}{2}\sum_{
        \scriptstyle x,y\in\RR
    \atop
        \scriptstyle x \neq y
    }c_{\RR,\lambda_\SS}(x,y)\left[f(y)-f(x)\right]^2,
    \qquad f: {\cal R} \rightarrow \mathds{R},
$$
where the {\it conductances} $c_{{\cal R}, \lambda_{\cal S}}(x,y)$
are given by
$$
    c_{\RR,\lambda_{\cal S}}(x,y)=\mu_{\RR}(x)w_{\RR,\lambda_\SS}(x,y)=\mu_{\RR}(y)w_{\RR,\lambda_\SS}(y,x).
$$
For $\lambda = + \infty$ we define $X_{{\cal R}, \lambda_{\cal S}}$
as the limiting process $X_{\cal R}$, rather than by removing all killing moves
from the limiting killed process $X^*_{{\cal R} \backslash {\cal S}}$.
In this case we will also write ${\cal D}_{\cal R}$ and $c_{\cal R}$
in place of ${\cal D}_{{\cal R}, \lambda_{\cal S}}$ and $c_{{\cal R}, \lambda_{\cal S}}$.
From Equation~\eqref{lodi}
as well as the continuity and monotonicity in $\lambda$
of the conductances $c_{{\cal R}, \lambda_{\cal S}}$,
which is inherited from that of the rates
$w_{{\cal R}, \lambda_{\cal S}}(x, y)$,
we get
\begin{prp}\label{pavia}
    When $\lambda$ goes to $+\infty$,
    $\gamma_{{\cal R}, \lambda_{\cal S}}$
    continuously decreases to $\gamma_{\cal R}$.
    In particular, it holds
    $$
        \gamma_{{\cal R}, \lambda_{\cal S}} \geq \gamma_{\cal R}
    $$
    for all $\lambda \geq 0$.
\end{prp}

We have from propositions~\ref{milano} and \ref{pavia}
that the ratio
$$
    \epsilon^*_{{\cal R}, \lambda_{\cal S}}
    = { 1 / \gamma_{{\cal R}, \lambda_{\cal S}} \over 1 / \phi^*_{{\cal R}, \lambda_{\cal S}}}
    = {\phi^*_{{\cal R}, \lambda_{\cal S}} \over \gamma_{{\cal R}, \lambda_{\cal S}}}
$$
is decreasing in $\lambda$.
When $\lambda = + \infty$ we will also write
\begin{equation}\label{bologna}
    \epsilon^*_{{\cal R}, {\cal S}}
    = {\phi^*_{{\cal R} \backslash {\cal S}} \over \gamma_{\cal R}}
\end{equation}
in place of $\epsilon^*_{{\cal R}, \lambda_{\cal S}}$
and we get,
for all $\lambda \geq 0$,
\begin{equation}\label{alessandria}
    \epsilon^*_{{\cal R}, \lambda_{\cal S}}
    \leq \epsilon^*_{{\cal R}, {\cal S}}.
\end{equation}
Hence $\epsilon^*_{{\cal R}, \lambda_{\cal S}}$
tends to be negligible under hypotheses $(H)$.

In this case all our soft measures are close together.
Let us introduce
some more notation to make this precise.
We denote by $h^*_{{\cal R}, \lambda_{\cal S}}$ the density
of $\mu^*_{{\cal R}, \lambda_{\cal S}}$ with respect to $\mu_{\cal R}$,
i.e.,
$$
    h^*_{{\cal R}, \lambda_{\cal S}}(x)
    = {\mu^*_{{\cal R}, \lambda_{\cal S}}(x) \over \mu_{\cal R}(x)}
    \,,
    \qquad x \in {\cal R}.
$$
In the special case $\lambda = +\infty$
we write $h^*_{{\cal R} \backslash {\cal S}}$
for $h^*_{{\cal R}, \lambda_{\cal S}}$:
$$
    h^*_{{\cal R} \backslash {\cal S}} :
    x \in {\cal R}
    \mapsto
    {\mu^*_{{\cal R} \backslash {\cal S}}(x) \over \mu_{\cal R}(x)}
    \,,
$$
which is supported on ${\cal R} \backslash {\cal S}$.
In particular, for any finite $\lambda$, we get
$$
    \phi^*_{{\cal R}, \lambda_{\cal S}}
    = {
        \langle
            h^*_{{\cal R}, \lambda_{\cal S}},
            -{\cal L}^*_{{\cal R}, \lambda_{\cal S}} h^*_{{\cal R}, \lambda_{\cal S}}
        \rangle_{\cal R}
    \over
        \left\|
            h^*_{{\cal R}, \lambda_{\cal S}}
        \right\|^2_{\cal R}
    }
    \,,
$$
and
$$
    \phi^*_{{\cal R} \backslash {\cal S}}
    = {
        \langle
            h^*_{{\cal R} \backslash {\cal S}},
            -{\cal L}^*_{{\cal R} \backslash {\cal S}} h^*_{{\cal R} \backslash {\cal S}}
        \rangle_{\cal R}
    \over
        \left\|
            h^*_{{\cal R} \backslash {\cal S}}
        \right\|^2_{\cal R}
    }
    \,.
$$
Since ${\cal L}^*_{{\cal R} \backslash {\cal S}}$
acts on functions of ${\cal R} \backslash {\cal S}$,
to make sense of
${\cal L}^*_{{\cal R} \backslash {\cal S}} h^*_{{\cal R} \backslash {\cal S}}$
we used in the previous equation
\begin{cvn}\label{rimini}
    Any function $f$ on ${\cal A} \subset {\cal X}$
    is identified with the function on its support ${\cal B} \subset {\cal A}$.
    In particular it also identified,
    for any ${\cal C} \supset {\cal B}$,
    with the function $g$ on ${\cal C}$
    that coincides with $f$ on ${\cal B}$
    and has zero value in ${\cal C} \backslash {\cal B}$.
\end{cvn}
\noindent
We can now compare our soft measures in $\ell^2(\mu_{\cal R})$
and $\ell^1(\mu_{\cal R})$ norms.
\begin{prp}\label{monza}
    For all $\lambda \in [0, +\infty]$,
    if $\epsilon_{\RR,\lambda_{\cal S}}^*<1$, then
    $$
        \var_{\mu_{\RR}}\lp h_{\RR,\lambda_\SS}^*\rp
        = \|h_{\RR,\lambda_\SS}^*-\mathds{1}_{\RR}\|^2_{\RR}
        \leq\frac{\epsilon_{\RR,\lambda_\SS}^*}{1-\epsilon_{\RR,\lambda_\SS}^*}
        \,.
    $$
    In particular, with $d_{\rm TV}$ the total variation distance,
    \begin{equation}\label{como}
        d_{\rm TV}\left(
            \mu^*_{{\cal R}, \lambda_{\cal S}}, \mu_{\cal R}
        \right)
        \leq {1 \over 2} \sqrt{\epsilon^*_{{\cal R}, \lambda_{\cal S}} \over 1 - \epsilon^*_{{\cal R}, \lambda_{\cal S}}}
        \,.
    \end{equation}
\end{prp}

\noindent{\it Proof:\/}
Using $h^*_{{\cal R}, \lambda_{\cal S}}$ as a test function
in the variational principle given by Equation~\eqref{lodi},
we get, for any finite $\lambda$,
\begin{align*}
    {\rm Var}_{\mu_{\cal R}}\bigl(h^*_{{\cal R}, \lambda_{\cal S}}\bigr)
    &\leq {1 \over \gamma_{{\cal R}, \lambda_{\cal S}}}
    {\cal D}_{{\cal R}, \lambda_{\cal S}}\bigl(h^*_{{\cal R}, \lambda_{\cal S}}\bigr)
    = {1 \over \gamma_{{\cal R}, \lambda_{\cal S}}}\langle
        h^*_{{\cal R}, \lambda_{\cal S}},
        -{\cal L}_{{\cal R}, \lambda_{\cal S}}h^*_{{\cal R}, \lambda_{\cal S}}
    \rangle_{\cal R} \\
    &= {1 \over \gamma_{{\cal R}, \lambda_{\cal S}}}\left(
        \langle
            h^*_{{\cal R}, \lambda_{\cal S}},
            -{\cal L}^*_{{\cal R}, \lambda_{\cal S}}h^*_{{\cal R}, \lambda_{\cal S}}
        \rangle_{\cal R}
        - \sum_{x \in {\cal R}} \mu_{\cal R}(x) e^*_{\cal R}(x) h^{*^2}_{{\cal R}, \lambda_{\cal S}}(x)
    \right) \\
    &\leq {1 \over \gamma_{{\cal R}, \lambda_{\cal S}}}
    \phi^*_{{\cal R}, \lambda_{\cal S}}
    \mu_{\cal R}\Bigr(h^{*^2}_{{\cal R}, \lambda_{\cal S}}\Bigl)
    = {\phi^*_{{\cal R}, \lambda_{\cal S}} \over \gamma_{{\cal R}, \lambda_{\cal S}}}\Bigl(
        1 + {\rm Var}_{\mu_{\cal R}}\bigr(h^*_{{\cal R}, \lambda_{\cal S}}\bigl)
    \Bigr).
\end{align*}
In the limiting case $\lambda = +\infty$
we obtain a similar estimate:
\begin{align*}
    {\rm Var}_{\mu_{\cal R}}\bigl(h^*_{{\cal R} \backslash {\cal S}}\bigr)
    &\leq {1 \over \gamma_{{\cal R}}}
    {\cal D}_{{\cal R}}\bigl(h^*_{{\cal R} \backslash {\cal S}}\bigr)\\
    &= {1 \over \gamma_{{\cal R}}}\left(
        {1 \over 2} \sum_{x, y \in {\cal R} \backslash {\cal S}} c_{\cal R}(x, y)\bigl[
            h^*_{{\cal R} \backslash {\cal S}}(x) - h^*_{{\cal R} \backslash {\cal S}}(y)
        \bigr]^2
        + \sum_{
            \scriptstyle x \in {\cal R} \backslash {\cal S}
            \atop
            \scriptstyle y \in {\cal S} \cap {\cal R}
        } c_{\cal R}(x, y) \bigl[
            h^*_{{\cal R} \backslash {\cal S}}(x)
        \bigr]^2
    \right)\\
    &\leq {1 \over \gamma_{{\cal R}}}\Biggl(
        \sum_{x \in {\cal R} \backslash {\cal S}} \mu_{\cal R}(x) h^*_{{\cal R} \backslash {\cal S}}(x)
        \sum_{y \in {\cal R} \backslash {\cal S}} w(x, y) \bigl[
            h^*_{{\cal R} \backslash {\cal S}}(x) - h^*_{{\cal R} \backslash {\cal S}}(y)
        \bigr] \\
        & \qquad\qquad + \sum_{x \in {\cal R} \backslash {\cal S}} \mu_{\cal R}(x) \sum_{y \in {\cal S}} w(x, y) \bigl[
            h^*_{{\cal R} \backslash {\cal S}}(x)
        \bigr]^2
    \Biggr)\\
    &= {1 \over \gamma_{{\cal R}}} \langle
        h^*_{{\cal R} \backslash {\cal S}},
        -{\cal L}^*_{{\cal R} \backslash {\cal S}}h^*_{{\cal R} \backslash {\cal S}}
    \rangle_{\cal R} \\
    &= {1 \over \gamma_{{\cal R}}}
    \phi^*_{{\cal R} \backslash {\cal S}}
    \mu_{\cal R}\Bigr(h^{*^2}_{{\cal R} \backslash {\cal S}}\Bigl)
    = {\phi^*_{{\cal R} \backslash {\cal S}} \over \gamma_{{\cal R}}}\Bigl(
        1 + {\rm Var}_{\mu_{\cal R}}\bigr(h^*_{{\cal R} \backslash {\cal S}}\bigl)
    \Bigr).
\end{align*}
Rearranging the terms, we get the desired inequality.
Equation~\eqref{como} follows then by Jensen inequality.
\qed

As a consequence,
as soon as $\lambda \gg \phi^*_{{\cal R}, \lambda_{\cal S}}$,
for example by setting $\lambda \gg \phi^*_{{\cal R} \backslash {\cal S}}$
(in virtue of~\eqref{neptuno_1}),
we obtain an asymptotic exponential law for $T_{\lambda_{\cal S}}$
when starting from $\mu^*_{{\cal R}, \lambda_{\cal S}}$ or $\mu_{\cl R}$.
If $\epsilon^*_{{\cal R}, \lambda_{\cal S}} \ll 1$,
then the random time
$\phi^*_{{\cal R}, \lambda_{\cal S}} T_{\lambda_{\cal S}}$
converges in law to an exponential variable or parameter~1.
\begin{prp}\label{gazza}
    For all $\lambda \geq \phi^*_{{\cal R} \backslash {\cal S}}$,
    or any $\lambda$ such that $\lambda \geq \phi^*_{{\cal R}, \lambda_{\cal S}}$,
    it holds
    $$
        e^{-t}
        \leq \mathds{P}_{\mu^*_{{\cal R}, \lambda_{\cal S}}}\left(
            T_{\lambda_{\cal S}}
            > {t \over \phi^*_{{\cal R}, \lambda_{\cal S}}}
        \right)
        \leq
        e^{-t} \left\{
            \exp\left(
                \sqrt{\phi^*_{{\cal R}, \lambda_{\cal S}} / \lambda}
            \right)
            + \exp\left(
                t - \sqrt{\lambda / \phi^*_{{\cal R}, \lambda_{\cal S}}}
            \right)
        \right\}
    $$
    as soon as $\sqrt{\phi^*_{{\cal R}, \lambda_{\cal S}} / \lambda} \leq t$,
    and, if $\epsilon^*_{{\cal R}, \lambda_{\cal S}} < 1$,
    $$
        \left|
            \mathds{P}_{\mu_{\cal R}}\left(
                T_{\lambda_{\cal S}}
                > {t \over \phi^*_{{\cal R}, \lambda_{\cal S}}}
            \right)
            - \mathds{P}_{\mu^*_{{\cal R}, \lambda_{\cal S}}}\left(
                T_{\lambda_{\cal S}}
                > {t \over \phi^*_{{\cal R}, \lambda_{\cal S}}}
            \right)
        \right|
        \leq
        {1 \over 2} \sqrt{
            \epsilon^*_{{\cal R}, \lambda_{\cal S}}
            \over 1 - \epsilon^*_{{\cal R}, \lambda_{\cal S}}
        }\,.
    $$
\end{prp}
\par\noindent{\bf Remark:}
As in \cite{BG} we write into brace parenthesis quantities
that go to one under suitable hypotheses.

\smallskip\par\noindent
\par\noindent
{\it Proof of the proposition:\/}
The first inequality is a consequence of the fact
that $T_{\lambda_{\cal S}} \geq T^{\cal R}_{\lambda_{\cal S}}$
and that, starting from $\mu^*_{{\cal R}, \lambda_{\cal S}}$,
the latter is an exponential variable
of parameter $\phi^*_{{\cal R}, \lambda_{\cal S}}$.
As far as the second one is concerned,
it holds, with $\ell_{\cal S}$ the local time in ${\cal S}$,
$$
    T_{\lambda_{\cal S}}
    \leq T^{\cal R}_{\lambda_{\cal S}}
    + \ell_{\cal S}\bigl(T_{\lambda_{\cal S}}\bigr)
$$
so that,
for any $0 \leq \beta \leq 1$,
\begin{align*}
    \mathds{P}_{\mu^*_{{\cal R}, \lambda_{\cal S}}}\left(
        T_{\lambda_{\cal S}} > {t \over \phi^*_{{\cal R}, \lambda_{\cal S}}}
    \right)
    &\leq \mathds{P}_{\mu^*_{{\cal R}, \lambda_{\cal S}}}\left(
        T^{\cal R}_{\lambda_{\cal S}} > {(1 - \beta) t \over \phi^*_{{\cal R}, \lambda_{\cal S}}}
    \right)
    + \mathds{P}_{\mu^*_{{\cal R}, \lambda_{\cal S}}}\left(
        \ell_{\cal S}\bigl(T_{\lambda_{\cal S}}\bigr)
        > {\beta t \over \phi^*_{{\cal R}, \lambda_{\cal S}}}
    \right)\\
    &= e^{- (1 - \beta) t} + e^{-\lambda \beta t / \phi^*_{{\cal R}, \lambda_{\cal S}}}.
\end{align*}
Choosing
$$
    \beta = {1 \over t} \sqrt{
        \phi^*_{{\cal R}, \lambda_{\cal S}} \over \lambda
    } \leq 1,
$$
this gives the desired result.
Finally,
by considering an optimal coupling
between two random variables $X(0)$ and $X'(0)$
of law $\mu_{\cal R}$ and $\mu^*_{{\cal R}, \lambda_{\cal S}}$,
such that $X(0) \neq X'(0)$ with probability
(recall Proposition~\ref{monza})
$$
    d_{\rm TV}\bigl(\mu_{\cal R}, \mu^*_{{\cal R}, \lambda_{\cal S}}\bigr)
    \leq {1 \over 2} \sqrt{
        \epsilon^*_{{\cal R}, \lambda_{\cal S}}
        \over 1 - \epsilon^*_{{\cal R}, \lambda_{\cal S}}
    }
$$
it holds
$$
    \mathds{P}_{\mu_{\cal R}}\left(
        \phi^*_{{\cal R}, \lambda_{\cal S}} T_{\lambda_{\cal S}} > t
    \right)
    \leq {1 \over 2} \sqrt{
        \epsilon^*_{{\cal R}, \lambda_{\cal S}}
        \over 1 - \epsilon^*_{{\cal R}, \lambda_{\cal S}}
    }
    + \mathds{P}_{\mu^*_{{\cal R}, \lambda_{\cal S}}}\left(
        \phi^*_{{\cal R}, \lambda_{\cal S}} T_{\lambda_{\cal S}} > t
    \right)
$$
and
$$
    \mathds{P}_{\mu^*_{{\cal R}, \lambda_{\cal S}}}\left(
        \phi^*_{{\cal R}, \lambda_{\cal S}} T_{\lambda_{\cal S}} > t
    \right)
    \leq {1 \over 2} \sqrt{
        \epsilon^*_{{\cal R}, \lambda_{\cal S}}
        \over 1 - \epsilon^*_{{\cal R}, \lambda_{\cal S}}
    }
    + \mathds{P}_{\mu_{\cal R}}\left(
        \phi^*_{{\cal R}, \lambda_{\cal S}} T_{\lambda_{\cal S}} > t
    \right).
$$
\qed

It follows that we also have an asymptotic exponential law
for $T_{\lambda_{\cal S}}$
when starting from any probability measure $\nu$
such that the total variation distance
between $\mu_{\cal R}$ and $\nu$ goes to zero.
An example of such a distribution $\nu$
is given by the law of $X(T_{\kappa_{\cal R}})$
if $\kappa \ll \gamma_{\cal R}/ \ln \chi_{\cal R}$
(recall Equation~\eqref{uovo}).
Hence, if a probability measure $\pi$ is such that
$\mathds{P}_\pi(T_{\kappa_{\cal R}} > T_{\lambda_{\cal S}}) \ll 1$,
then $T_{\lambda_{\cal S}}$ has an asymptotic exponential law
of parameter $\phi^*_{{\cal R}, \lambda_{\cal S}}$.
\begin{prp}\label{tortora}
    For all $x \in {\cal X}$ and all $\kappa > 0$
    it holds
    \begin{equation}\label{pfff}
        d_{\rm TV}\bigl(
            \mu_{\cal R},
            \mathds{P}_x \left(
                X(T_{\kappa_{\cal R}}) = \cdot
            \right)
        \bigr)
        \leq {\kappa \over \gamma_{\cal R}} \left(
            1 + \left[
                \ln {
                    \gamma_{\cal R} \sqrt{\chi_{\cal R}}
                    \over 2\kappa
                }
            \right]_+
        \right).
    \end{equation}
    Also, for any probability $\pi$ on ${\cal X}$, by setting
    $$
        \delta = \mathds{P}_\pi\left(
            T_{\lambda_{\cal S}} < T_{\kappa_{\cal R}}
        \right)
    $$
    and
    $$
        \eta
        = \delta
        + {\kappa \over \gamma_{\cal R}} \left(
            1 + \left[
                \ln{
                    \gamma_{\cal R} \sqrt{\chi_{\cal R}}
                    \over 2\kappa
                }
            \right]_+
        \right)
        + {1 \over 2}\sqrt{
            \epsilon^*_{{\cal R}, \lambda_{\cal S}}
            \over
            1 - \epsilon^*_{{\cal R}, \lambda_{\cal S}}
        }\,,
    $$
    it holds, for all $\lambda \geq \phi^*_{{\cal R} \backslash {\cal S}}$,
    or any $\lambda \geq \phi^*_{{\cal R}, \lambda_{\cal S}}$,
    $$
        e^{-t}\Bigl\{
            1 - \eta e^t
        \Bigr\}
        \leq \mathds{P}_\pi\left(
            T_{\lambda_{\cal S}} > {t \over \phi^*_{{\cal R}, \lambda_{\cal S}}}
        \right)
        \leq e^{-t} \left\{
            \exp\left(
                2 \sqrt{\phi^*_{{\cal R}, \lambda_{\cal S}} \over \kappa \wedge \lambda}
            \right)
            + \left[
                \eta
                + 2 \exp\left(
                    - \sqrt{\kappa \wedge \lambda \over \phi^*_{{\cal R}, \lambda_{\cal S}}}
                \right)
            \right] e^t
        \right\}.
    $$
\end{prp}

\par\noindent
{\it Proof:\/} Let us denote by $X^{\cal R}$ the trace of $X$
on ${\cal R}$.
Note that we already gave to this process an heavier notation:
$$
    X^{\cal R} = X_{{\cal R}, 0_{\cal S}}.
$$
Since $X(T_{\kappa_{\cal R}}) \in {\cal R}$, the law
of $X(T_{\kappa_{\cal R}})$
is that of $X^{\cal R}(T^{\cal R}_{\kappa_{\cal R}})$
with
$$
    T^{\cal R}_{\kappa_{\cal R}}
    = \ell_{\cal R}\bigl(T_{\kappa_{\cal R}}\bigr).
$$
Since the relaxation time of $X^{\cal R}$ is smaller than
or equal to that of $X_{\cal R}$ (recall Proposition~\ref{pavia}),
by using reversibility
to estimate the $\ell^2(\mu_{\cal R})$ distance
between (the densities with respect to $\mu_{\cal R}$ of)
the law of $X^{\cal R}(s)$
and $\mu_{\cal R}$,
then Jensen inequality to compare
the $\ell^1$ and $\ell^2$ distances,
we have
$$
    d_{\rm TV}\left(
        \mathds{P}_x\Bigl(X^{\cal R}(s) = \cdot\Bigr),
        \mu_{\cal R}
    \right)
    \leq {\sqrt{\chi_{\cal R}} \over 2} e^{-\gamma_{\cal R}s}.
$$
Since $T^{\cal R}_{\kappa_{\cal R}}$ is an exponential variable
of parameter $\kappa$, it then holds, for all non-negative $s$,
$$
    d_{\rm TV}\left(
        \mathds{P}_x\Bigl(X(T_{\kappa_{\cal R}}) = \cdot\Bigr),
        \mu_{\cal R}
    \right)
    \leq \kappa s + {\sqrt{\chi_{\cal R}} \over 2} e^{-\gamma_{\cal R}s}.
$$
Optimizing in $s \geq 0$, we get \eqref{pfff}.

Next, denoting by $\nu$ the law of $X(T_{\kappa_{\cal R}})$
with starting distribution $\pi$,
it holds, for all $t \geq 0$,
$$
    \mathds{P}_\nu\left(
        T_{\lambda_{\cal S}} > {t \over \phi^*_{{\cal R}, \lambda_{\cal S}}}
    \right)
    \leq \mathds{P}_\pi\left(
        T_{\lambda_{\cal S}} < T_{\kappa_{\cal R}}
    \right)
    + \mathds{P}_\pi\left(
        T_{\lambda_{\cal S}} > {t \over \phi^*_{{\cal R}, \lambda_{\cal S}}}
    \right).
$$
By using Proposition~\ref{gazza}
it follows that
$$
    \mathds{P}_\pi\left(
        T_{\lambda_{\cal S}} > {t \over \phi^*_{{\cal R}, \lambda_{\cal S}}}
    \right)
    \geq e^{-t} - {1 \over 2}\sqrt{
        \epsilon^*_{{\cal R}, \lambda_{\cal S}}
        \over
        1 - \epsilon^*_{{\cal R}, \lambda_{\cal S}}
    } - {\kappa \over \gamma_{\cal R}} \left(
        1 + \left[
            \ln {
                \gamma_{\cal R} \sqrt{\chi_{\cal R}}
                \over 2\kappa
            }
        \right]_+
    \right) - \delta.
$$
Finally, for all $\theta > 0$, we get
$$
    \mathds{P}_\pi\left(
        T_{\lambda_{\cal S}} > {t \over \phi^*_{{\cal R}, \lambda_{\cal S}}}
    \right)
    \leq
    \mathds{P}_\pi\left(
        T_{\lambda_{\cal S}} < T_{\kappa_{\cal R}}
    \right)
    + \mathds{P}_\pi\left(
        T_{\lambda_{\cal S}} \wedge T_{\kappa_{\cal R}} > \theta
    \right)
    + \mathds{P}_\nu\left(
        T_{\lambda_{\cal S}}
        > {t \over \phi^*_{{\cal R}, \lambda_{\cal S}}} - \theta
    \right).
$$
By using Proposition~\ref{gazza}
and the fact that $T_{\kappa_{\cal R}} \wedge T_{\lambda_{\cal S}}$
is dominated by an exponential variable
of parameter $\kappa \wedge \lambda$
we get
\begin{align*}
    \mathds{P}_\pi\left(
        T_{\lambda_{\cal S}} > {t \over \phi^*_{{\cal R}, \lambda_{\cal S}}}
    \right)
    &\leq \delta
    + e^{-(\kappa \wedge \lambda) \theta}
    + {1 \over 2} \sqrt{
        \epsilon^*_{{\cal R}, \lambda_{\cal S}}
        \over 1 - \epsilon^*_{{\cal R}, \lambda_{\cal S}}
    }
    + {\kappa \over \gamma_{\cal R}} \left(
        1 + \left[
            \ln {
                \gamma_{\cal R} \sqrt{\chi_{\cal R}}
                \over 2\kappa
            }
        \right]_+
    \right)
    + \exp\left[
        - \sqrt{\lambda / \phi^*_{{\cal R}, \lambda_{\cal S}}}
    \right]\\
    &\qquad + e^{-t} \exp\left[
        \phi^*_{{\cal R}, \lambda_{\cal S}} \theta
        + \sqrt{\phi^*_{{\cal R}, \lambda_{\cal S}} / \lambda}
    \right],
\end{align*}
which gives our last desired inequality
by choosing
$\theta = (\phi^*_{{\cal R}, \lambda_{\cal S}} (\kappa \wedge \lambda))^{- 1 / 2}$.
\qed

\subsection{Soft capacities}\label{livorno}
To define our soft capacities
we extend the conductance network $({\cal X}, c)$,
with conductances
$$
    c(x, y) = c(y, x) = \mu(x) w(x, y),
    \qquad x, y \in {\cal X}, \qquad x \neq y,
$$
which is associated with the process $X$,
into a larger network $(\tilde {\cal X}, \tilde c)$
by attaching dangling edges $(r, \bar r)$
and $(s, \breve s)$ to each $r \in {\cal R}$
and $s \in {\cal S}$.
Given $\kappa$ and $\lambda$ we set to
$$
    \tilde c(r, \bar r) = \kappa \mu(r)
    \qquad{\rm and}\qquad
    \tilde c(s, \breve s) = \lambda \mu(s)
$$
the conductances of these new edges.
We call $\bar {\cal R}$ and $\breve {\cal S}$
the collections of these extra nodes:
$$
    \bar {\cal R} = \bigl\{
        \bar r : r \in {\cal R}
    \bigr\},
    \qquad
    \breve {\cal S} = \bigl\{
        \breve s : s \in {\cal S}
    \bigr\}.
$$
For any positive $\tilde k$ and $\tilde \lambda$,
by extending the probability measure $\mu$ defined on ${\cal X}$
into a measure $\tilde \mu$ defined on
$$
    \tilde {\cal X} = {\cal X} \cup \bar {\cal R} \cup \breve {\cal S}
$$
by
$$
    \tilde \mu(\tilde x)
    = \left\{
        \begin{array}{ll}
            \mu(x) & \hbox{if $\tilde x = x$ for some $x \in {\cal X}$,}\\
            \kappa \mu(r) / \tilde \kappa & \hbox{if $\tilde x = \bar r$ for some $r \in \bar {\cal R}$,}\\
            \lambda \mu(s) / \tilde \lambda & \hbox{if $\tilde x = \breve s$ for some $s \in \breve {\cal S}$,}
        \end{array}
    \right.
$$
the extended network $(\tilde X, \tilde c)$
is then associated with the process $\tilde X$
with rates $\tilde w(\tilde x, \tilde y)$ given,
for distinct $\tilde x$ and $\tilde y$ in $\tilde {\cal X}$,
by
$$
    \tilde w(\tilde x, \tilde y)
    = \left\{
        \begin{array}{ll}
            w(x, y) & \hbox{%
                if $\tilde x = x$ and $\tilde y = y$ for some $x$ and $y$ in ${\cal X}$,
            }\\
            \kappa & \hbox{%
                if $\tilde x = r$ and $\tilde y = \bar r$ for some $r$ in ${\cal R}$,
            }\\
            \tilde \kappa & \hbox{%
                if $\tilde x = \bar r$ and $\tilde y = r$ for some $r$ in ${\cal R}$,
            }\\
            \lambda & \hbox{%
                if $\tilde x = s$ and $\tilde y = \breve s$ for some $s$ in ${\cal S}$,
            }\\
            \tilde \lambda & \hbox{%
                if $\tilde x = \breve s$ and $\tilde y = s$ for some $s$ in ${\cal S}$,
            }\\
            0 & \hbox{%
                otherwise,
            }
        \end{array}
    \right.
$$
and which is reversible with respect to $\tilde \mu$.
(Note that $\tilde \mu$ is not a probability measure.)

Our soft capacities are then the $(\kappa, \lambda)$-capacities
of \cite{BG}.
\begin{dfn}\label{varese}
    The $(\kappa,\lambda)$-capacity $C_{\kappa}^{\lambda}\lp\RR,\SS\rp$
    is the capacity between the sets $\bar{\RR}$ and $\breve{\SS}$
    in the network $\lp\tilde\XX,\tilde{c}\rp$.
    According to Dirichlet principle, it is given by
    \begin{align}
        C_{\kappa}^{\lambda}\lp\RR,\SS\rp
        &=\min\left\{
            \frac{1}{2}\sum_{\tilde{x},\tilde{y}\in\tilde\XX}\tilde{c}\bigl(\tilde{x},\tilde{y}\bigl) \Bigl(
                \tilde{f}(\tilde{x})-\tilde{f}(\tilde{y})
            \Bigr)^2
            : \tilde{f}|_{\bar\RR}=1, \tilde{f}|_{\breve\SS}=0
        \right\} \nonumber\\
        &=\min_{f:\XX\to\R}\left\{
            \DD(f)
            + \kappa\sum_{r\in\RR}\mu(r)\bigl(f(r)-1\bigr)^2
            + \lambda\sum_{s\in\SS}\mu(s)\bigl(f(s)-0\bigr)^2
        \right\} \nonumber\\
        &=\min_{f:\XX\to\R}\left\{
            \DD(f)
            + \kappa\mu(\RR) {\rm E}_{\mu_{\RR}}\left[
                \bigl(f_{|_\RR}-1\bigr)^2
            \right]
            + \lambda\mu(\SS) {\rm E}_{\mu_{\SS}}\left[
                \bigl(f_{|_\SS}-0\bigr)^2
            \right]
        \right\}. \label{sondrio}
    \end{align}
\end{dfn}

\noindent {\bf Remarks:}
\begin{itemize}
\item[i.] Since this definition depends
    on the conductances only,
    while it is naturally associated
    with the Markov process $\tilde X$,
    it does not depends on the choice
    of $\tilde \kappa$ and $\tilde \lambda$.
    $C_\kappa^\lambda({\cal R}, {\cal S})$
    depends on two parameters only,
    $\kappa$ and $\lambda$,
    in which it is increasing.
\item[ii.] Our soft capacities satisfy
    a two-sided variational principle.
    On one side they are given by Definition~\ref{varese}
    as the infimum of some functional,
    and any test function
    will provide an upper bound.
    On the other side they are, by Thomson principle,
    the supremum of another functional on unitary flows
    from $\bar{\RR}$ to $\breve{\SS}$,
    which are antisymmetric functions
    of oriented edges with null divergence in $\XX$
    and total divergence in $\bar {\cal R}$
    and $\breve {\cal S}$ equal to 1 and -1 respectively,
    i.e., on antisymmetric functions
    $\tilde\psi : \tilde\XX \times \tilde\XX \rightarrow \RR$
    such that $\tilde \psi(\tilde x, \tilde y) = 0$
    as soon as $\tilde c(\tilde x, \tilde y) = 0$
    and
    $$ \text{div}_x\tilde\psi=\sum_{\tilde{x}\in\tilde\XX}\tilde\psi(x,\tilde{x})=0,
        \qquad x \in {\cal X},
    $$
    while
    $$
        \sum_{r \in {\cal R}} {\rm div}_{\bar r} \tilde \psi
        = \sum_{r \in {\cal R}} \tilde \psi \bigl(\bar r, r\bigr)
        = 1
        = - \sum_{s \in {\cal S}} {\rm div}_{\breve s} \tilde \psi
        = \sum_{s \in {\cal S}} \tilde \psi(s, \breve s).
    $$
    With
    $$
        \widetilde\DD\lp
            \tilde\psi
        \rp = \frac{1}{2}\sum_{\tilde{x},\tilde{y}\in\tilde\XX}\frac{\tilde\psi(\tilde{x},\tilde{y})^2}{\tilde{c}(\tilde{x},\tilde{y})}
    $$
    the energy dissipated in $\bigl(\tilde{\cal X}, \tilde c\bigr)$
    by the flow $\tilde\psi$ and with $\widetilde\Psi_1\bigl(\bar\RR,\breve\SS\bigr)$
    the set of unitary flows from $\bar\RR$ to $\breve\SS$, we get
    \begin{equation}\label{bergamo}
        C_\kappa^\lambda\left(
            \RR,\SS
        \right) = \max_{\tilde\psi\in\widetilde\Psi_1(\bar\RR,\breve\SS)} \widetilde\DD\lp\tilde\psi\rp^{-1}.
    \end{equation}
    Any test flow provides then a lower bound on $C_\kappa^\lambda\lp\RR,\SS\rp$.
\item[iii.] By referring to $\tilde X$,
    the stopping time $T_{\lambda_{\cal S}}$
    can also be defined as the hitting time of $\breve {\cal S}$.
    By defining $T_{\kappa_{\cal R}}$ in a symmetric way,
    the minimum in \eqref{sondrio}
    is realized by the equilibrium potential $V_\kappa^\lambda$ given by
    $$
        V_\kappa^\lambda(x) = \mathds{P}_x\left(
                T_{\kappa_{\cal R}} < T_{\lambda_{\cal S}}
            \right),
        \qquad x \in {\cal X},
    $$
    whereas the maximum in \eqref{bergamo} is realized by the associated normalized current.
\end{itemize}

As an application of Dirichlet principle
we obtain an upper bound on our soft capacities,
or of the asymptotic given for the spectral $\gamma$
in Theorem~\ref{albero}.
Let us define
$$
    \phi_\kappa^\lambda
    = {C_\kappa^\lambda({\cal R}, {\cal S}) \over \mu({\cal R}) \mu({\cal S})}
$$
and recall Equation~\eqref{bologna}.
\begin{lmm}\label{ravenna}
    If $\mu({\cal S}) \geq \mu({\cal R})$
    and $\epsilon^*_{{\cal R},{\cal S}} \leq 1$,
    then, for all $\kappa \geq 0$,
    it holds
    \begin{equation}\label{cuneo}
        \phi_\kappa^\lambda
        \leq  \mu({\cal S})^{-1} \phi^*_{{\cal R} \backslash {\cal S}}\left\{
            1 + \kappa  / \gamma_{\cal R}
        \over
            1 - \epsilon^*_{{\cal R}, {\cal S}}
        \right\}
        \qquad \hbox{and} \qquad
        \mu(S)^{-1} \leq 2.
    \end{equation}
In particular, under the hypotheses of Theorem~\ref{albero}
it holds $\phi_\kappa^\lambda = O\bigl(\phi^*_{{\cal R} \backslash {\cal S}}\bigr)$.
\end{lmm}

\par\noindent
{\it Proof:\/}
By taking $h^*_{\RR\backslash\SS}$ as test function in Equation \eqref{sondrio},
and following Convention~\ref{rimini}, we get
\begin{equation*}
    C_{\kappa}^\lambda(\RR,\SS)
    \leq \DD\bigl(h^*_{\RR\backslash\SS}\bigr) + \kappa\mu(\RR)\var_{\mu_\RR}\bigl(h^*_{\RR\backslash\SS}\bigr).
\end{equation*}
The second term is straightforwardly bounded by using Proposition \ref{monza}
which reads
$$
    \kappa\mu(\RR)\var_{\mu_\RR}\bigl(h^*_{\RR\backslash\SS}\bigr)
    \leq \mu({\cal R}) \phi^*_{{\cal R} \backslash {\cal S}}
    {\kappa / \gamma_{\cal R} \over 1 - \epsilon^*_{{\cal R}, {\cal S}}}
    \,.
$$
 As far as the first term is concerned we write
\begin{align*}
    \DD\lp h^*_{\RR\backslash\SS}\rp&=\frac{1}{2}\sum_{x,y\in\RR\backslash\SS}c(x,y)\bigl[h^*_{\RR\backslash\SS}(x)-h^*_{\RR\backslash\SS}(y)\bigr]^2+\sum_{\substack{x\in\RR\backslash\SS,\\ y\in\SS}}c(x,y)h^{*^2}_{\RR\backslash\SS}(x)\\
    &=\mu(\RR)\langle h^*_{\RR\backslash\SS}, -\LL^*_{\RR\backslash\SS}h^*_{\RR\backslash\SS}\rangle_\RR
    =\mu(\RR)\phi^*_{\RR\backslash\SS}\bigl\|h^*_{\RR\backslash\SS}\bigr\|^2_\RR
    =\mu(\RR)\phi^*_{\RR\backslash\SS}\lp1+\var_{\mu_\RR}\bigl( h^*_{\RR\backslash\SS}\bigr)\rp\\
    &\leq\mu(\RR)\phi^*_{\RR\backslash\SS}
    \lp1+\frac{\epsilon^*_{\RR,\SS}}{1-\epsilon^*_{\RR,\SS}}\rp
    =\mu({\cal R}) \phi^*_{{\cal R}\backslash {\cal S}} {1 \over 1 - \epsilon^*_{{\cal R}, {\cal S}}}\,,
\end{align*}
where the inequality is again an application Proposition \ref{monza}.
Observing that ${\cal X} = {\cal R} \cup {\cal S}$
and $\mu({\cal S}) \geq \mu({\cal R})$
imply $\mu({\cal S}) \geq 1/ 2$,
this leads to~\eqref{cuneo}.
The last comment is only based on the fact
that hypotheses $(H)$ include $\epsilon^*_{{\cal R}, {\cal S}} \ll 1$
and $\kappa \ll \gamma_{\cal R}$ is assumed in Theorem~\ref{albero}.
\qed

\begin{rmk}\label{asti}
    All the estimates of Section~\ref{sole}
    are ``symmetric'' in ${\cal R}$ and ${\cal S}$.
    This means that the same results hold
    when the role of ${\cal R}$ and ${\cal S}$
    are exchanged as well as those of $\kappa$ and $\lambda$.
\end{rmk}
\noindent
We note however that there is an asymmetry
in hypotheses $(H)$. It will start
to play a role in the next section.

\section{Proof of Theorem~\ref{albero}} \label{luna}
\subsection{Spectral gap estimates}
The statement contained in Theorem \ref{albero} concerning the spectral gap $\gamma$
is a consequence of the quantitative Theorem \ref{piacenza}
together with Lemma \ref{ravenna}, estimates \eqref{neptuno_1} and \eqref{alessandria} as well as Remark~\ref{asti}.
Note also that hypotheses $(H)$ include
$\phi^*_{{\cal R} \backslash {\cal S}} \ll \gamma_{\cal S}$.
This is where the asymmetry of $(H)$ matters.

\begin{thm}\label{piacenza}
    The spectral gap $\gamma$ satisfies
    \begin{equation}\label{forli}
        \gamma\geq\phi_{\kappa}^{\lambda}
        \lb1+\frac{\kappa+\phi_{\kappa}^{\lambda}\lp1+\mu\lp\RR\cap\SS\rp\rp}{\gamma_{\RR}}+\frac{\lambda+\phi_{\kappa}^{\lambda}
        \lp1+\mu\lp\RR\cap\SS\rp\rp}{\gamma_{\SS}}\rb^{-1}.
    \end{equation}
    In addition, when $\epsilon^*_{{\cal R}, \lambda_{\cal S}}$
    and $\epsilon^*_{{\cal S}, \kappa_{\cal R}}$ are smaller than $1$,
    we obtain the converse inequality
    \begin{equation}\label{verbania}
        \gamma\leq \phi_{\kappa}^{\lambda}\lb
            1 - \frac{\phi_{\RR,\lambda_\SS}^*}{\kappa}-\frac{\phi_{\SS,\kappa_\RR}^*}{\lambda}
            - \frac{1}{2}\sqrt{\frac{\epsilon_{\RR,\lambda_\SS}^*}{1-\epsilon_{\RR,\lambda_\SS}^*}}
            -\frac{1}{2}\sqrt{\frac{\epsilon_{\SS,\kappa_\RR}^*}{1-\epsilon_{\SS,\kappa_\RR}^*}}
        \rb^{-2}
    \end{equation}
    as soon as the braced sum is positive.
\end{thm}

\par\noindent
{\it Proof of the lower bound:\/}
We recall that $\gamma$ satisfies the variational principle
\begin{equation}
\gamma:=\min_{\substack{f:\XX\to\R, \\ \var_{\mu}(f)\neq 0}}\frac{\DD(f)}{\var_{\mu}(f)}\,,\label{modena}
\end{equation}
We look for a Poincaré inequality, i.e. an inequality of the form
\begin{equation*}
    \forall f\in\ell^2(\mu), \quad \var_{\mu}(f)\leq k\DD\lp f\rp,
\end{equation*}
which is equivalent, by mean of Equation \eqref{modena}, to a lower bound on the spectral gap.
Let $Z, Z'$ two i.i.d. random variables with same law $\mu$.
The variance of $f$ can be written as
\begin{equation}
    \var_{\mu}(f)
    =\text{E}_{\mu}\left[ f^2\right]-\text{E}_{\mu}^2\left[ f\right]
    =\frac{1}{2}\text{E}\left[ \lp f(Z)-f(Z')\rp^2\right].\label{var1}
\end{equation}
From the r.h.s. of (\ref{var1}) we get
\allowdisplaybreaks
\begin{align}
\var_{\mu}(f)
&\leq\frac{1}{2}\text{E}\left[\lp f(Z)-f(Z')\rp^2\mathds{1}_{\left\{ Z,Z'\in\RR\right\}}\right] +
\frac{1}{2}\text{E}\left[\lp f(Z)-f(Z')\rp^2\mathds{1}_{\left\{Z,Z'\in\SS\right\}}\right]\nonumber\\
&\qquad +\text{E}\left[\lp f(Z)-f(Z')\rp^2\mathds{1}_{\left\{Z\in \RR, Z'\in\SS\right\}}\right]\nonumber\\
&=\mu^2(\RR)\var_{\mu_{\RR}}\lp f_{|\RR}\rp+\mu^2(\SS)
\var_{\mu_{\SS}}\lp f_{|\SS}\rp\nonumber\\
&\qquad +\mu(\RR)\mu(\SS)\lp \text{E}_{\mu_{\RR}}\left[f_{|\RR}^2\right]
 + \text{E}_{\mu_{\SS}}\left[f_{|\SS}^2\right]-2\text{E}_{\mu_{\RR}}\left[f_{|\RR}\right]\text{E}_{\mu_{\SS}}\left[f_{|\SS}\right]\rp\nonumber\\
&=\mu(\RR)\Bigl(\mu(\RR)+\mu(\SS)\Bigr)\var_{\mu_{\RR}}\lp f_{|\RR}\rp + \mu(\SS)\Bigl(\mu(\RR)+\mu(\SS)\Bigr)\var_{\mu_{\SS}}\lp f_{|\SS}\rp\nonumber\\
&\qquad +\mu(\RR)\mu(\SS)\Bigl(
    \text{E}^2_{\mu_{\RR}}\left[f_{|\RR}\right]+\text{E}^2_{\mu_{\SS}}\left[f_{|\SS}\right]
    -2\text{E}_{\mu_{\RR}}\left[f_{|\RR}\right]\text{E}_{\mu_{\SS}}\left[f_{|\SS}\right]
\Bigr) \nonumber\\
&=\mu(\RR)\Bigl(\mu(\RR)+\mu(\SS)\Bigr)\var_{\mu_{\RR}}\lp f_{|\RR}\rp + \mu(\SS)\Bigl(\mu(\RR)+\mu(\SS)\Bigr)\var_{\mu_{\SS}}\lp f_{|\SS}\rp\nonumber\\
&\qquad +\mu(\RR)\mu(\SS)\Bigl(\text{E}_{\mu_{\RR}}\left[f_{|\RR}\right]-\text{E}_{\mu_{\SS}}\left[f_{|\SS}\right]\Bigr)^2.\label{diffsqa}
\end{align}
As a test function in Dirichlet principle \eqref{sondrio} we pick
\begin{equation}
    \tilde{f}=\frac{f-\text{E}_{\mu_{\SS}}\left[f_{|\SS}\right]}{\text{E}_{\mu_{\RR}}\left[f_{|\RR}\right]-\text{E}_{\mu_{\SS}}\left[f_{|\SS}\right]}\label{testf},
\end{equation}
which is such that $\text{E}_{\mu_{\RR}}\left[\tilde{f}_{|{\cal R}}\right]=1$
and $\text{E}_{\mu_{\SS}}\left[\tilde{f}_{|{\cal S}}\right]=0$. By plugging (\ref{testf}) in (\ref{sondrio}) we get
\begin{align*}
C_{\kappa}^{\lambda}(\RR,\SS)&\leq\DD(\tilde{f})+\kappa\mu(\RR)\var_{\mu_{\RR}}\lp\tilde{f}_{|\RR}\rp+\lambda\mu(\SS)\var_{\mu_{\SS}}\lp\tilde{f}_{|\SS}\rp\nonumber\\
&=\Bigl(
    \text{E}_{\mu_{\RR}}\left[f_{|\RR}\right]-\text{E}_{\mu_{\SS}}\left[f_{|\SS}\right]
\Bigr)^{-2}\Bigl(
    \DD(f)+\kappa\mu(\RR)\var_{\mu_{\RR}}\lp f_{|\RR}\rp+\lambda\mu(\SS)\var_{\mu_{\SS}}\lp f_{|\SS}\rp
\Bigr)
\end{align*}
which actually is an upper bound on $\lp \text{E}_{\mu_{\RR}}\left[f_{|\RR}\right]-\text{E}_{\mu_{\SS}}\left[f_{|\SS}\right]\rp^{2}$.
Bounding from above the last term in (\ref{diffsqa})
we then get
\allowdisplaybreaks
\begin{align*}
\var_{\mu}(f)&\leq \mu(\RR)\Bigl(\mu(\RR)+\mu(\SS)\Bigr)\var_{\mu_{\RR}}\lp f_{|\RR}\rp+\mu(\SS)\Bigl(\mu(\RR)+\mu(\SS)\Bigr)\var_{\mu_{\SS}}\lp f_{|\SS}\rp\\
&+\frac{\mu(\RR)\mu(\SS)}{C_{\kappa}^{\lambda}\lp\RR,\SS\rp}\Bigl(\DD(f)+\kappa\mu(\RR)\var_{\mu_{\RR}}(f_{|{\RR}})+\lambda\mu(\SS)\var_{\mu_{\SS}}\lp f_{|\SS}\rp\Bigr)\\
&=\frac{\DD(f)}{\phi_{\kappa}^{\lambda}}+\mu(\RR)\var_{\mu_{\RR}}(f_{|{\RR}})\Biggl(\frac{\kappa}{\phi_{\kappa}^{\lambda}} +\lp 1 +\mu\lp\RR\cap \SS\rp\rp\Biggr)\\
& \qquad + \mu(\SS)\var_{\mu_{\SS}}(f_{|{\SS}})\Biggl(\frac{\lambda}{\phi_{\kappa}^{\lambda}}+\lp 1+\mu\lp\RR\cap \SS\rp\rp\Biggr)\\
&\leq \frac{\DD(f)}{\phi_{\kappa}^{\lambda}}+\frac{\mu(\RR)\DD_{\RR}(f_{|{\RR}})}{\gamma_{\RR}}\lp \frac{\kappa}{\phi_{\kappa}^{\lambda}}+\lp 1+\mu\lp\RR\cap \SS\rp\rp\rp+\frac{\mu(\SS)\DD_{\SS}(f_{|{\SS}})}{\gamma_{\SS}}\lp \frac{\lambda}{\phi_{\kappa}^{\lambda}}+\lp 1+\mu\lp\RR\cap \SS\rp\rp\rp\nonumber\\
&\leq\frac{\DD(f)}{\phi_{\kappa}^{\lambda}}+\frac{\DD(f)}{\gamma_{\RR}}\lp \frac{\kappa}{\phi_{\kappa}^{\lambda}}+\lp 1+\mu\lp\RR\cap \SS\rp\rp\rp+\frac{\DD(f)}{\gamma_{\SS}}\lp \frac{\lambda}{\phi_{\kappa}^{\lambda}}+\lp 1+\mu\lp\RR\cap \SS\rp\rp\rp
\end{align*}
where, in the second inequality,
the variational characterization of spectral gaps $\gamma_\RR$ and $\gamma_\SS$ given by Equation \eqref{reggio_emilia}  has been used,
whereas in the third one we exploited the following fact:
\begin{equation*}
\DD(f)\geq\mu(\RR)\DD_{\RR}(f_{|{\RR}}), \qquad \DD(f)\geq\mu(\SS)\DD_{\SS}(f_{|{\SS}}).
\end{equation*}
Factorising ${\cal D}(f) / \phi_\kappa^\lambda$ we get the Poincaré inequality
\begin{equation*}
    \var_{\mu}(f)
    \leq \frac{\DD(f)}{\phi_{\kappa}^{\lambda}}\lp
        1+\frac{\kappa+\phi_{\kappa}^{\lambda} \lp1+\mu\lp\RR\cap\SS\rp\rp}{\gamma_{\RR}}
        +\frac{\lambda+\phi_{\kappa}^{\lambda} \lp1+\mu\lp\RR\cap\SS\rp\rp}{\gamma_{\SS}}
    \rp.
\end{equation*}
\qed

\smallskip\par
Before giving the proof of the upper bound,
we observe that the proof of this Poincaré inequality carries over
to the case of any finite covering
of irreducible sets
${\cal X} = \cup_{i < m} {\cal R}_i$.
\begin{prp}\label{giggiola}
    If ${\cal R}_0$, ${\cal R}_1$, \dots, ${\cal R}_{m - 1}$
    are subsets of ${\cal X}$ such that $X_{{\cal R}_i}$
    is irreducible for each $i < m$,
    and $\kappa_0$, $\kappa_1$, \dots, $\kappa_{m - 1}$
    are any positive numbers
    then, with
    $$
        \phi_i^j = \phi_j^i = {
            C_{\kappa_i}^{\kappa_j}({\cal R}_i, {\cal R}_j)
        \over
            \mu({\cal R}_i) \mu({\cal R}_j)
        }\,,
        \qquad i \neq j, \qquad i, j < m,
    $$
    $$
        {1 \over \phi}
        = {1 \over 2} \sum_{i, j : i \neq j} {1 \over \phi_i^j}
    $$
    and, for each $i < m$,
    $$
        {1 \over \phi_i}
        = \sum_{j \neq i} {1 \over \phi_i^j}
        \leq {1 \over \phi},
    $$
    it holds
    $$
        \gamma
        \geq \phi \left\{
            1 + \sum_{i < m} {\phi \over \gamma_i} \left(
                {\kappa_i \over \phi_i}
                + \sum_{j < m} \mu({\cal R}_j)
            \right)
        \right\}^{-1}
        \geq \phi \left\{
            1 + \sum_{i < m} {
                \kappa_i  + \phi \sum_{j < m} \mu({\cal R}_j)
            \over
                \gamma_i
            }
        \right\}^{-1}.
    $$
\end{prp}

\par\noindent
{\it Proof of the upper bound in Theorem~\ref{piacenza}:\/}
We want to exploit the variational principle \eqref{modena}, so we look for an upper bound on $\DD(f)$ and a lower bound on $\var_\mu(f)$
for a suitable function $f$.
From \eqref{var1} we get
\allowdisplaybreaks
\begin{align}
\var_{\mu}(f)&\geq \frac{1}{2}\text{E}\left[
    \bigl( f(Z)-f(Z')\bigr)^2\lp\mathds{1}_{\left\{Z\in \RR, Z'\in \SS\right\}}+\mathds{1}_{\left\{Z\in \SS, Z'\in \RR\right\}}-\mathds{1}_{\left\{Z,Z'\in \RR\cap\SS\right\}}\rp
\right]\nonumber\\
&=\text{E}\left[ \bigl( f(Z)-f(Z')\bigr)^2\mathds{1}_{\left\{Z\in \RR, Z'\in \SS\right\}}\right]
-\frac{1}{2}\text{E}\left[ \bigl( f(Z)-f(Z')\bigr)^2\mathds{1}_{\left\{Z,Z'\in \RR\cap\SS\right\}}\right]\nonumber\\
&=\mu(\RR)\mu(\SS)\lp \text{E}_{\mu_{\RR}}\left[ f_{|{\RR}}^2\right]+\text{E}_{\mu_{\SS}}\left[ f_{|{\SS}}^2\right]\rp-2\mu(\RR)\mu(\SS)\text{E}_{\mu_{\RR}}\left[f_{|{\RR}}\right]\text{E}_{\mu_{\SS}}\left[ f_{|{\SS}}\right]\nonumber\\
&\qquad -\mu^2(\RR\cap \SS)\var_{\mu_{\RR\cap \SS}}\lp f_{|{\RR\cap \SS}}\rp\nonumber\\ &=\mu(\RR)\mu(\SS)\Bigl(\var_{\mu_{\RR}}\lp f_{|{\RR}}\rp+\text{E}_{\mu_{\RR}}^2\left[ f_{|{\RR}}\right]\Bigr)+\mu(\RR)\mu(\SS)\Bigl(\var_{\mu_{\SS}}\lp f_{|{\SS}}\rp+\text{E}_{\mu_{\SS}}^2\left[ f_{|{\SS}}\right]\Bigr)\nonumber\\
&\qquad -2\mu(\RR)\mu(\SS)\text{E}_{\mu_{\RR}}\left[ f_{|{\RR}}\right]\text{E}_{\mu_{\SS}}\left[ f_{|{\SS}}\right]-\mu^2(\RR\cap \SS)\var_{\mu_{\RR\cap \SS}}\lp f_{|{\RR\cap\SS}}\rp\nonumber\\
&=\mu(\RR)\mu(\SS)\Bigl(\var_{\mu_{\RR}}\lp f_{|{\RR}}\rp+\var_{\mu_{\SS}}\lp f_{|{\SS}}\rp\Bigr)-\mu^2(\RR\cap \SS)\var_{\mu_{\RR\cap \SS}}\lp f_{|{\RR\cap \SS}}\rp\label{var2}\\
&\qquad +\mu(\RR)\mu(\SS)\Bigl( \text{E}_{\mu_{\RR}}\left[ f_{|{\RR}}\right]-\text{E}_{\mu_{\SS}}\left[ f_{|{\SS}}\right]\Bigr)^2\nonumber
\end{align}
To get a suitable lower bound of \eqref{var2}, we write
\begin{align}
\mu^2(\RR\cap \SS)\var_{\mu_{\RR\cap \SS}}\lp f_{|\RR\cap \SS}\rp
&=\frac{1}{2}\text{E}\left[ \bigl( f(Z)-f(Z')\bigr)^2\mathds{1}_{\left\{Z,Z'\in \RR\cap\SS\right\}}\right]\nonumber\\
&\leq\frac{1}{2}\text{E}\left[ \bigl( f(Z)-f(Z')\bigr)^2\mathds{1}_{\left\{Z,Z'\in \RR\right\}}\right]\nonumber\\
&=\mu^2(\RR)\var_{\mu_{\RR}}\lp f_{|{\RR}}\rp
\leq \mu(\RR)\mu(\SS)\var_{\mu_{\RR}}\lp f_{|{\RR}}\rp,\label{var3}
\end{align}
assuming, without loss of generality, $\mu(\RR)\leq\mu(\SS)$.
Then, by using \eqref{var3} to bound the second term in \eqref{var2},
dropping one term and taking as test function the equilibrium potential $V_\kappa^\lambda$, we get
\begin{equation}
\var_{\mu}\bigl( V_{\kappa}^{\lambda}\bigr)\geq \mu(\RR)\mu(\SS)\Bigl(     \text{E}_{\mu_{\RR}}\left[ V_{\kappa}^{\lambda}|_{\RR}\right]-\text{E}_{\mu_{\SS}}\left[ V_{\kappa}^{\lambda}|_{\SS}\right]\Bigr)^2.\label{varclev}
\end{equation}

\goodbreak
Next, we have the two following simples lemmas.
\begin{lmm}\label{lemmino1}
If $\epsilon^*_{{\cal R}, \lambda_{\cal S}} < 1$, then
\begin{equation}
\text{E}_{\mu_{\RR}} \ls V_{\kappa}^{\lambda}|_\RR\rs\geq 1-\frac{\phi_{\RR,\lambda_\SS}^*}{\kappa}-\frac{1}{2}\sqrt{\frac{\epsilon_{\RR,\lambda_\SS}^*}{1-\epsilon_{\RR,\lambda_\SS}^*}}\nonumber.
\end{equation}
\end{lmm}

\par\noindent
{\it Proof:\/}
 Denoting by $\sigma_\kappa$ an exponential time of rate $\kappa$ that is independent of $X$, it holds
\begin{equation}
\text{E}_{\mu_{\RR}}\ls V_{\kappa}^{\lambda}|_{\RR}\rs=\P_{\mu_{\RR}}\lp T_{\kappa_{\RR}}<T_{\lambda_{\SS}}\rp=\P_{\mu_{\RR}}\lp \ell_{\RR}\lp T_{\kappa_{\RR}}\rp<\ell_{\RR}\lp T_{\lambda_{\SS}}\rp\rp=
\P_{\mu_{\RR}}\lp\sigma_{\kappa}< T_{\lambda_\SS}^{\RR}\rp.\nonumber
\end{equation}
Let $\P^{opt}$ be an optimal coupling at time $t=0$
between $X^{\mu_\RR}$ and $X^{\mu_{\RR,\lambda_{\cal S}}^*}$
that start with initial distribution $\mu_{\cal R}$
and $\mu^*_{{\cal R}, \lambda_{\cal S}}$
and evolve jointly if they have the same starting configuration.
By the union bound, Proposition \eqref{monza}
 and the convexity of $x\mapsto\frac{1}{1+x}$ we get
\begin{align}
\P_{\mu_{\RR}}\lp\sigma_{\kappa}<T_{\lambda_\SS}^{\cal R}\rp
&\geq\P_{\mu_{\RR,\lambda}^*}\lp\sigma_{\kappa}<T_{\lambda_\SS}^\RR\rp-\P^{opt}\lp X^{\mu_\RR}(0)\neq X^{\mu^*_{\RR,\lambda_{\cal S}}}(0)\rp\nonumber\\
&\geq\frac{\kappa}{\kappa+\phi_{\RR,\lambda_\SS}^*}-\frac{1}{2}\sqrt{\frac{\epsilon_{\RR,\lambda_\SS}^*}{1-\epsilon_{\RR,\lambda_\SS}^*}}\geq 1-\frac{\phi_{\RR,\lambda_\SS}^*}{\kappa}-\frac{1}{2}\sqrt{\frac{\epsilon_{\RR,\lambda_\SS}^*}{1-\epsilon_{\RR,\lambda_\SS}^*}}\nonumber
\end{align}
\qed

\begin{lmm}\label{lemmino2}
If $\epsilon^*_{{\cal S}, \kappa_{\cal R}} < 1$, then
\begin{equation}
E_{\mu_{\SS}}\left[ V_{\kappa}^{\lambda}|_{\SS}\right]\leq \frac{1}{2}\sqrt{\frac{\epsilon_{S,\kappa_\RR}^*}{1-\epsilon_{S,\kappa_\RR}^*}} +\frac{\phi_{S,\kappa_\RR}^*}{\lambda}\nonumber.
\end{equation}
\end{lmm}
\par\noindent
The proof of this lemma is similar to that of the previous one
and we omit it.
By plugging the two bounds provided by Lemma \ref{lemmino1} and Lemma \ref{lemmino2} into \eqref{varclev} we obtain
\begin{equation}\label{arezzo}
    \var_{\mu}\lp V_{\kappa}^{\lambda}\rp\geq\mu(\RR)\mu(\SS)\lb1-\frac{\phi_{\RR,\lambda_\SS}^*}{\kappa}-\frac{\phi_{\SS,\kappa_\RR}^*}{\lambda}-\frac{1}{2}\sqrt{\frac{\epsilon_{\RR,\lambda_\SS}^*}{1-\epsilon_{\RR,\lambda_\SS}^*}}-\frac{1}{2}\sqrt{\frac{\epsilon_{\SS,\kappa_\RR}^*}{1-\epsilon_{\SS,\kappa_\RR}^*}}\rb^2
\end{equation}
as soon as the braced quantity is positive.
The upper bound on $\DD(V_\kappa^\lambda)$ is straightforward:
since the functional associated with the Dirichlet principle,
which reaches its minimum $C_\kappa^\lambda({\cal R}, {\cal S})$
in $V_\kappa^\lambda$, is larger than ${\cal D}$,
we get
\begin{equation}\label{grosseto}
    \DD\bigl( V_{\kappa}^{\lambda}\bigr)\leq C_{\kappa}^{\lambda}\lp \RR,\SS\rp.
\end{equation}
Inequality~\eqref{verbania} follows from Equation~\eqref{modena}
together with \eqref{arezzo} and \eqref{grosseto}.
\qed

\subsection{Exit rate estimates}\label{brune}

The prove the last part of Theorem~\ref{albero}
we first give sharp upper and lower bounds
for $\phi^*_{\cl R, \lambda_{\cl S}}$.
These are sharp bounds in the sense
that they have the same asymptotic behaviour
as soon as $\kappa$ and $\lambda$ are chosen
in the right window.

\begin{thm}\label{pisa}
The exit rate $\phi_{\RR,\lambda_\SS}^*$ satisfies
\begin{equation*}
    \phi_{\RR,\lambda_\SS}^*
    \geq {C_\kappa^\lambda({\cal R}, {\cal S}) \over \mu({\cal R})} \left\{
        1 - 2\mu({\cal S})^{-1} \phi^*_{{\cal R}, \lambda_{\cal S}} / \lambda
    \over
        \bigl(1 - \phi^*_{{\cal R}, \lambda_{\cal S}} / \lambda\bigr)^2
    \right\}
    \lb1+\frac{\kappa+\phi_{\kappa}^{\lambda}\lp1+\mu\lp\RR\cap\SS\rp\rp}{\gamma_{\RR}}+\frac{\lambda+\phi_{\kappa}^{\lambda}\lp1+\mu\lp\RR\cap\SS\rp\rp}{\gamma_{\SS}}\rb^{-1}\!.
\end{equation*}
In addition, when $\epsilon^*_{{\cal R}, \lambda_{\cal S}} < 1$,
we get the converse inequality
\begin{equation*}
    \phi_{\RR,\lambda_\SS}^*
    \leq \frac{C_{\kappa}^{\lambda}\lp\RR,\SS\rp}{\mu(\RR)}
    \lb 1+\frac{\phi_{\RR,\lambda_\SS}^*}{\lambda}\rb
    \lb1-\frac{\phi_{\RR,\lambda_\SS}^*}{\kappa}-\frac{1}{2}\sqrt{\frac{\epsilon_{\RR,\lambda_\SS}^*}{1-\epsilon_{\RR,\lambda_\SS}^*}}\rb^{-2}
\end{equation*}
as soon as the last braced sum is positive.
\end{thm}

\par\noindent
{\it Proof of the lower bound:\/}
We consider the process $\tilde X$ introduced in Section~\ref{livorno}
in the limiting case $\tilde \kappa = + \infty$.
In this case $\tilde {\cal X}$ can be redefined
as $\tilde {\cal X} = {\cal X} \cup \breve {\cal S}$
and we call $\tilde {\cal L}$
the generator of this process $\tilde X$
on $\tilde {\cal X}$.
According with the notation of Section~\ref{novara},
we denote by $\tilde\phi_{\XX,\lambda_\SS}^*$
and $\tilde\mu_{\XX,\lambda_\SS}^*$ respectively
the smallest eigenvalue and the corresponding normalized left eigenvector
(i.e., the quasi-stationary distribution associated with $\XX$ and $\lambda_{\cal S}$)
of $-\tilde\LL^*_{\XX,\lambda_\SS}$,
which is defined by
$$
    \bigl(\tilde{\cal L}^*_{{\cal X}, \lambda_{\cal S}}f\bigr)(x)
    = - \lambda_{\cal S}(x) f(x) + \sum_{y \in {\cal X}} w(x, y)[f(y) - f(x)],
    \qquad f: {\cal X} \rightarrow \mathds{R},
    \qquad x \in {\cal X}.
$$

For any probability $\nu$ on ${\cal X}$,
$\tilde \phi^*_{{\cal X}, \lambda_{\cal S}}$
is the decay rate of the survival probability
$\mathds{P}_\nu\bigl(T_{\lambda_{\cal S}} > t\bigr)$:
\begin{equation*}
    \lim_{t\to\infty}-\frac{1}{t}\ln\P_\nu\lp
         T_{\lambda_\SS} > t
    \rp
    = \tilde\phi^*_{\XX,\lambda_\SS}.
\end{equation*}
Since
\begin{equation*}
    \P_{\mu_{\RR,\lambda_\SS}^*} \lp
        T_{\lambda_\SS} > t
    \rp
    \geq \P_{\mu_{\RR,\lambda_\SS}^*}\Bigl(
        \ell_{\RR}(T_{\lambda_\SS}) > t
    \Bigr)
    = \P_{\mu_{\RR,\lambda_\SS}^*}\Bigl(
        T_{\lambda_\SS}^\RR >t
    \Bigr)
    = e^{-t\phi_{\RR,\lambda_\SS}^*},
\end{equation*}
it follows that
\begin{equation}
    \phi_{\RR,\lambda_\SS}^*\geq\tilde\phi_{\XX,\lambda_\SS}^*.\label{ivrea}
\end{equation}

To give a lower bound on $\tilde\phi^*_{{\cal X}, \lambda_{\cal S}}$,
we use its variational representation
\begin{equation}\label{vercelli}
    \tilde\phi_{\XX,\lambda_\SS}^*
    =\min_{f:\XX\to \R}\frac{\bigl\langle
        f, -\tilde {\cal L}^*_{{\cal X}, \lambda_{\cal S}} f
    \bigr\rangle}
    {\| f\|^2}
    =\min_{f:\XX\to \R}\frac{
        {\cal D}(f) + \|\sqrt{\lambda_{\cal S}} f\|^2
    }
    {\| f\|^2}
    \,.
\end{equation}
Indeed, for any $f : {\cal X} \rightarrow \mathds{R}$,
\begin{equation*}
    \left\langle
         f,-\tilde\LL_{\XX,\lambda_\SS}^*f
    \right\rangle
    = \frac{1}{2}\sum_{x,y\in\XX}\mu(x)w(x,y)\bigl[f(x)-f(y)\bigr]^2
    + \sum_{x\in\XX}\mu(x)\lambda_\SS(x)f^2(x).
\end{equation*}
Since the minimum is reached in
$$
    \tilde{h}_{\XX,\lambda_\SS}^*:
    x \in {\cal X}
    \mapsto \tilde\mu^*_{\XX,\lambda_\SS}(x) / \mu(x).
$$
it holds
\begin{equation}
\tilde\phi_{\XX,\lambda_\SS}^*\geq\frac{\DD\bigl(\tilde{h}_{\XX,\lambda_\SS}^*\bigr)}{\bigl\|\tilde{h}_{\XX,\lambda_\SS}^*\bigr\|^2}=\frac{\bigl\|\tilde{h}_{\XX,\lambda_\SS}^*\bigr\|^2-1}{\bigl\|\tilde{h}_{\XX,\lambda_\SS}^*\bigl\|^2}\frac{\DD\lp\tilde{h}_{\XX,\lambda_\SS}^*\rp}{\var_{\mu}\lp\tilde{h}_{\XX,\lambda_\SS}^*\rp}\geq \lp1-\frac{1}{\bigl\|\tilde{h}_{\XX,\lambda_\SS}^*\bigr\|^2}\rp\gamma.\label{testnorm3}
\end{equation}

Theorem~\ref{piacenza}
provides a lower bound on $\gamma$
and it remains to give a lower bound
on $\|\tilde h^*_{{\cal X}, \lambda_{\cal S}}\|$.
To estimate this norm
we restrict the sum to ${\cal X} \backslash {\cal S}$
and we use Jensen's inequality:
\begin{align}
    \bigl\|\tilde{h}_{\XX,\lambda_\SS}^*\bigr\|^2
    & \geq \mu({\cal X} \backslash {\cal S}) \sum_{x\in {\cal X} \backslash {\cal S}} \mu_{{\cal X} \backslash {\cal S}}(x)\left(
        \frac{\tilde\mu_{\XX,\lambda_\SS}^*(x)}{\mu(x)}
    \right)^2  \nonumber \\
    & \geq \mu({\cal X} \backslash {\cal S})\left(
        \sum_{x\in {\cal X} \backslash {\cal S}} \mu_{{\cal X} \backslash {\cal S}}(x)
        \frac{\tilde\mu_{\XX,\lambda_\SS}^*(x)}{\mu(x)}
    \right)^2
    =\frac{\bigl(\tilde\mu_{\XX,\lambda_\SS}^*({\cal X} \backslash {\cal S})\bigr)^2}{\mu({\cal X} \backslash {\cal S})}
    \,.\label{testnorm}
\end{align}
Now, the last identity of Lemma~\ref{saturno} reads in this case
$\tilde{\phi}_{\XX,\lambda_\SS}^*=\lambda\tilde\mu_{\XX,\lambda_\SS}^*(\SS)$.
Therefore, using Inequality~\eqref{ivrea} again,
\begin{equation}\label{testnorm2}
\tilde\mu^*_{\XX,\lambda_\SS}({\cal X} \backslash {\cal S})
= 1-\frac{\tilde\phi_{\XX,\lambda_\SS}^*}{\lambda}
\geq 1-\frac{\phi_{\RR,\lambda_\SS}^*}{\lambda}
\,.
\end{equation}
It follows from  \eqref{testnorm} and \eqref{testnorm2} that
\begin{align}
    1 - {1 \over \bigl\|\tilde h^*_{{\cal X}, \lambda_{\cal S}}\bigr\|^2}
    &\geq 1 - {
        1 - \mu(\SS)
    \over
        \bigl(1 - \phi^*_{{\cal R}, \lambda_{\cal S}} / \lambda\bigr)^2
    }
    = {
        \mu({\cal S})
        - 2  \phi^*_{{\cal R}, \lambda_{\cal S}} / \lambda
        + \bigl(\phi^*_{{\cal R}, \lambda_{\cal S}} / \lambda\bigr)^2
    \over
        \bigl(1 - \phi^*_{{\cal R}, \lambda_{\cal S}} / \lambda\bigr)^2
    }\nonumber\\
    &\geq \mu({\cal S}) {
        1 - 2  \mu({\cal S})^{-1} \phi^*_{{\cal R}, \lambda_{\cal S}} / \lambda
    \over
        \bigl(1 - \phi^*_{{\cal R}, \lambda_{\cal S}} / \lambda\bigr)^2
    }.\label{bari}
\end{align}
Finally from \eqref{ivrea}, \eqref{testnorm3} and \eqref{bari} together with \eqref{forli} we get
the desired result.
\qed

\smallskip\par\noindent
{\it Proof of the upper bound:\/}
The proof is made of two steps: first we estimate $\tilde\phi_{\XX,\lambda_\SS}^*$
via a variational principle,
then we look for an upper bound of $\phi_{\RR,\lambda_\SS}^*$ in terms of $\tilde\phi_{\XX,\lambda_\SS}^*$.
By taking as test function in \eqref{vercelli} the equilibrium potential $V_\kappa^\lambda$ we get
\begin{equation}
\tilde\phi_{\XX,\lambda_\SS}^*\leq \frac{\DD\lp V_\kappa^\lambda\rp+ \bigl\|\sqrt{\lambda_\SS} V_\kappa^\lambda\bigr\|^2}{\bigl\|V_\kappa^\lambda\bigr\|^2}\leq\frac{C_{\kappa}^{\lambda}\lp\RR,\SS\rp}{\bigl\|V_{\kappa}^{\lambda}\bigr\|^2}.\label{abba}
\end{equation}
The denominator in \eqref{abba} is lower bounded by
\allowdisplaybreaks
\begin{align}\label{abba2}
    \bigl\|V_{\kappa}^{\lambda}\bigr\|^2\geq\mu\lp\RR\rp\sum_{x\in\RR}\mu_{\RR}(x)V_{\kappa}^{\lambda}(x)^2\geq\mu\lp\RR\rp\mu_{\RR}\bigl(V_{\kappa}^{\lambda}|_\RR\bigr)^2\geq \mu\lp\RR\rp\lb1-\frac{\phi_{\RR,\lambda_\SS}^*}{\kappa}-\frac{1}{2}\sqrt{\frac{\epsilon_{\RR,\lambda_\SS}^*}{1-\epsilon_{\RR,\lambda_\SS}^*}}\rb^2,
\end{align}
where we used Jensen inequality, Lemma \ref{lemmino1}
and the assumption that the braced sum was positive.

To bound $\phi^*_{{\cal R}, \lambda_{\cal S}}$
with $\tilde \phi^*_{{\cal X}, \lambda_{\cal S}}$,
we write, for any $0 \leq \beta \leq 1$
and denoting again by $\sigma_\lambda$ an
exponential time of rate $\lambda$
that is independent of $X$,
\begin{align*}
    \P_{\mu_{\RR,\lambda_\SS}^*}\bigl( T_{\lambda_\SS}> t\bigr)
    &\leq \P_{\mu_{\RR,\lambda}^*}\Bigl(
        T_{\lambda_\SS}^\RR+\sigma_{\lambda}
        >\beta t+(1-\beta)t
    \Bigr)\\
    &\leq \P_{\mu_{\RR,\lambda}^*}\Bigl(
        T_{\lambda_\SS}^\RR > \beta t
    \Bigr) + \P_{\mu_{\RR,\lambda}^*}\bigl(
        \sigma_{\lambda} > (1-\beta)t
    \bigr)\\
    &= \exp\bigl(-\phi_{\RR,\lambda_{\cal S}}^*\beta t\bigr) + \exp\bigl(-\lambda(1-\beta)t\bigr).
\end{align*}
Recall that $\tilde \phi^*_{{\cal X}, \lambda_{\cal S}}$
is the decay rate of the survival probability
$\P_{\mu_{\RR,\lambda_\SS}^*}\bigl( T_{\lambda_\SS}> t\bigr)$.
If we choose $\beta$ such that
$\phi_{\RR,\lambda_{\cal S}}^*\beta\leq\lambda(1-\beta)$,
then $\phi^*_{{\cal R}, \lambda_{\cal S}}\beta$
is the decay rate of the right-hand side in the previous inequality
and we can conclude
$\tilde\phi^*_{{\cal X}, \lambda_{\cal S}} \geq \beta \phi^*_{{\cal R}, \lambda_{\cal S}}$.
We then choose
$\beta={\lambda} / ({\phi_{\RR,\lambda_{\cal S}}^*+\lambda})$, which gives us
\begin{equation}
    \phi_{\RR,\lambda_{\cal S}}^*\leq\frac{\tilde\phi_{\XX, \lambda_{\cal S}}^*}{\beta}=\tilde\phi_{\XX, \lambda_{\cal S}}^*
    \lp 1+\frac{\phi_{\RR,\lambda_{\cal S}}^*}{\lambda}\rp.\label{conco}
\end{equation}
The thesis follows from \eqref{conco}, \eqref{abba} and \eqref{abba2}.
\qed

\smallskip\par
By Lemma~\ref{ravenna},
estimates \eqref{neptuno_1} and \eqref{alessandria}
as well as Remark~\ref{asti},
we only have to prove that
$\phi^*_{\cl R, \lambda_{\cl S}} \ds E_{\mu_{\cl R}}\bigl[
    T_{\lambda_{\cl S}}
\bigr]$ goes to one
to conclude the proof of Theorem~\ref{albero}.
Since
$$
    \ds E_{\mu_{\cl R}}\bigr[
        T_{\lambda_{\cl S}}
    \bigl]
    = \ds E_{\mu_{\cl R}}\bigr[
        \ell_{\cl R}(T_{\lambda_{\cl S}})
        + \ell_{\cl S \setminus \cl R}(T_{\lambda_{\cl S}})
    \bigl]
    \leq \ds E_{\mu_{\cl R}}\bigr[
        \ell_{\cl R}(T_{\lambda_{\cl S}})
    \bigr]
    + \ds E_{\mu_{\cl R}}\bigr[
        \ell_{\cl S}(T_{\lambda_{\cl S}})
    \bigl]
    = \ds E_{\mu_{\cl R}}\bigr[
        T^{\cl R}_{\lambda_{\cl S}}
    \bigr]
    + {1 \over \lambda}\,,
$$
we get
$$
    \phi^*_{\cl R, \lambda_{\cl S}}
    \ds E_{\mu_{\cl R}}\bigr[
        T_{\lambda_{\cl S}}
    \bigl]
    \leq 1 + {\phi^*_{\cl R, \lambda_{\cl S}} \over \lambda}
$$
from estimate~\eqref{neptuno_3}.
This upper bound goes to one,
as a consequence of estimate~\eqref{neptuno_1}.

To get a lower bound,
we write
$$
    \phi^*_{\cl R, \lambda_{\cl S}}
    \ds E_{\mu_{\cl R}}\bigr[
        T_{\lambda_{\cl S}}
    \bigl]
    = \int_0^{+\infty} \ds P_{\mu_{\cl R}}\bigl(
        \phi^*_{\cl R, \lambda_{\cl S}} T_{\lambda_{\cl S}}
        > t
    \bigr) \, dt
$$
and recall Proposition~\ref{gazza}:
the integrand goes to $e^{-t}$.
Fatou's lemma gives the desired lower bound
and this concludes the proof.
\qed

\section{Proof of Theorem~\ref{gheppio}}\label{pranzo}
We will use the stopping time
$$
    \tau = T_{\kappa_{\cal R}} \wedge T_{\lambda_{\cal S}}
$$
to build $T^*$.
The laws of $X(T_{\kappa_{\cal R}})$
and $X(T_{\lambda_{\cal S}})$
are indeed close to $\mu_{\cal R}$ and $\mu_{\cal S}$.
But conditioning
on $\{\tau = T_{\kappa_{\cal R}}\}$
or $\{\tau = T_{\lambda_{\cal S}}\}$
introduces correlations that can make the law
of $X(\tau)$ in general delicate to control.
The law of $X(T_{\kappa_{\cal R}})$
can be written as a convex combination
of those of $X(\tau)$ conditioned on $\{\tau = T_{\kappa_{\cal R}}\}$
and $X(T_{\kappa_{\cal R}})$
conditioned on $\{T_{\kappa_{\cal R} > T_{\lambda_{\cal S}}}\}$:
it holds, for all $x$ in ${\cal X}$,
and with
$$
    \alpha(x) = \mathds{P}_x\left(
        T_{\kappa_{\cal R}} < T_{\lambda_{\cal S}}
    \right),
$$
$$
    \mathds{P}_x\left(
        X(T_{\kappa_{\cal R}}) = \cdot
    \right)
    = \alpha(x) \mathds{P}_x\left(
        X(\tau) = \cdot
        \bigm| \tau = T_{\kappa_{\cal R}}
    \right)
    + (1 - \alpha(x)) \mathds{P}_x\left(
        X(T_{\kappa_{\cal R}}) = \cdot
        \bigm| T_{\kappa_{\cal R}} > T_{\lambda_{\cal S}}
    \right).
$$
By Proposition~\ref{tortora}
the first and last of these three distributions
are close to $\mu_{\cal R}$
for a small enough $\kappa$.
We can then estimate the total variation distance
between $\mu_{\cal R}$ and the law of $X(\tau)$
conditioned on $\{\tau = T_{\kappa_{\cal R}}\}$:
$$
    d_{\rm TV}\Bigl(
        \mu_{\cal R},
        \mathds{P}_x\left(
            X(\tau) = \cdot
            \bigm| \tau = T_{\kappa_{\cal R}}
        \right)
    \Bigr)
    \leq {2 - \alpha(x) \over \alpha(x)}
    {\kappa \over \gamma_{\cal R}} \left(
        1 + \left[
            \ln {
                \gamma_{\cal R} \sqrt{\chi_{\cal R}}
                \over 2\kappa
            }
        \right]_+
    \right),
$$
so that, if $\alpha(x) \geq 1 / 2$,
$$
    d_{\rm TV}\Bigl(
        \mu_{\cal R},
        \mathds{P}_x\left(
            X(\tau) = \cdot
            \bigm| \tau = T_{\kappa_{\cal R}}
        \right)
    \Bigr)
    \leq {3\kappa \over \gamma_{\cal R}} \left(
        1 + \left[
            \ln {
                \gamma_{\cal R} \sqrt{\chi_{\cal R}}
                \over 2\kappa
            }
        \right]_+
    \right).
$$

This suggests the following construction.
We set $\tau^0 = 0$
and, for all $i \geq 0$, we define
$$
    \tau^{i + 1}_{\cal R} = \tau^i + T_{\kappa_{\cal R}} \circ \Theta_{\tau^i},
    \qquad
    \tau^{i + 1}_{\cal S} = \tau^i + T_{\lambda_{\cal S}} \circ \Theta_{\tau^i},
    \qquad
    \tau^{i + 1} = \tau^i_{\cal R} \wedge \tau^i_{\cal S}.
$$
We call $I$ the first index $i > 0$ such that
either
$$
    \tau^i = \tau^i_{\cal R}
    \quad\hbox{and}\quad
    \alpha(X(\tau^{i - 1})) \geq {1 \over 2}
$$
or
$$
    \tau^i = \tau^i_{\cal S}
    \quad\hbox{and}\quad
    \alpha(X(\tau^{i - 1})) \leq {1 \over 2}.
$$
We finally define
$$
    T^*_{\cal R} = \tau^I_{\cal R},
    \qquad
    T^*_{\cal S} = \tau^I_{\cal S},
    \qquad
    T^* = T^*_{\cal R} \wedge T^*_{\cal S}.
$$

The random variable $I$ is stochastically dominated
by a geometric random variable of mean~2,
each random variable $(\tau^i - \tau^{i -1})$ is stochastically dominated
by an exponential random variable of rate $\kappa \wedge \lambda$
and it holds
$$
    T^* = \sum_{i > 0} (\tau^i - \tau^{i - 1}) \mathds{1}_{\{I \geq i\}},
$$
so that, observing that for $(\tau^i - \tau^{i - 1})$ and $\{I \geq i\}$
are independent for all $i > 0$,
$$
    \mathds{E}_x\left[
        T^*
    \right] \leq {2 \over \kappa \wedge \lambda}\,,
$$
for all $x$ in ${\cal X}$.
By construction it holds
$$
    d_{\rm TV}\Bigl(
        \mu_{\cal R},
        \mathds{P}_x\left(
            X(\tau) = \cdot
            \bigm| T^* = T^*_{\cal R}
        \right)
    \Bigr)
    \leq {3\kappa \over \gamma_{\cal R}} \left(
        1 + \left[
            \ln {
                \gamma_{\cal R} \sqrt{\chi_{\cal R}}
                \over 2\kappa
            }
        \right]_+
    \right).
$$
With
$$
    \eta
    = {3\kappa \over \gamma_{\cal R}} \left(
        1 + \left[
            \ln{
                \gamma_{\cal R} \sqrt{\chi_{\cal R}}
                \over 2\kappa
            }
        \right]_+
    \right)
    + {1 \over 2}\sqrt{
        \epsilon^*_{{\cal R}, \lambda_{\cal S}}
        \over
        1 - \epsilon^*_{{\cal R}, \lambda_{\cal S}}
    }\,,
$$
the triangular inequality gives
$$
    d_{TV}\Bigl(
        \mu^*_{{\cal R}, \lambda_{\cal S}},
        \mathds{P}_x\left(
            X(T^*) = \cdot
            \bigm| T^* = T^*_{\cal R}
        \right)
    \Bigr) \leq \eta\,.
$$
Using Proposition~\ref{gazza},
for all $\lambda \geq \phi^*_{{\cal R} \backslash {\cal S}}$,
or any $\lambda \geq \phi^*_{{\cal R}, \lambda_{\cal S}}$,
we get
$$
    \mathds{P}_x\left(
        T_{\lambda_{\cal S}} \circ \Theta_{T^*} > {t \over \phi^*_{{\cal R}, \lambda_{\cal S}}}
        \biggm| T^* = T^*_{\cal R}
    \right)
    \geq e^{-t}\Bigl\{
        1 - \eta e^t
    \Bigr\}
$$
and
$$
    \mathds{P}_x\left(
        T_{\lambda_{\cal S}} \circ \Theta_{T^*} > {t \over \phi^*_{{\cal R}, \lambda_{\cal S}}}
        \biggm| T^* = T^*_{\cal R}
    \right)
    \leq e^{-t} \left\{
        \exp\left(
            \sqrt{\phi^*_{{\cal R}, \lambda_{\cal S}} \over \lambda}
        \right)
        + \left[
            \eta
            + \exp\left(
                - \sqrt{\lambda \over \phi^*_{{\cal R}, \lambda_{\cal S}}}
            \right)
        \right] e^t
    \right\}.
$$
Finally, if
$$
    \delta = \mathds{P}_\nu\left(
        T_{\lambda_{\cal S}} < T_{\kappa_{\cal R}}
    \right),
$$
then, by definition of $\alpha$,
$$
    {1 \over 2}\,\mathds{P}_\nu\left(
        \alpha(X(0)) \leq {1 \over 2}
    \right)
    \leq \mathds{P}_\nu\left(
        T_{\lambda_{\cal S}} < T_{\kappa_{\cal R}}
        \Bigm| \alpha(X(0)) \leq {1 \over 2}
    \right) \mathds{P}_\nu\left(
        \alpha(X(0)) \leq {1 \over 2}
    \right)
    \leq \delta
$$
and
$$
   \mathds{P}_\nu\left(
        T^* = T^*_{\cal R} = T_{\kappa_{\cal R}} < T_{\lambda_{\cal S}}
    \right)
    \geq  1 - \mathds{P}_\nu\left(
        \alpha(X(0)) \leq {1 \over 2}
        \hbox{ or } T_{\lambda_{\cal S}} < T_{\kappa_{\cal R}}
    \right)
    \geq 1 - (2 \delta + \delta) = 1 - 3 \delta.
$$
Since the symmetrical statements hold when exchanging
the roles of ${\cal R}$ and ${\cal S}$
as well as those of $\kappa$ and $\lambda$,
this concludes the proof.
\qed

\section{Proof of Theorem~\ref{corteo}}\label{telefono}
By Proposition~\ref{tortora}
and Remark~\ref{asti},
the law of $X(T_1) = X(T_{\lambda_{\cal S}})$
is close to $\mu_{\cal S}$,
as the law of $X(T_2)$ is close to $\mu_{\cal R}$.
Hence, it will be sufficient to find $\theta$
such that our time averages on time scale $\theta$
are close to the expected values computed
with $\mu_{\cal R}$ before $T_1$.
Theorem~\ref{corteo} follows then from the quantitative
\begin{prp}
    If $\epsilon^*_{{\cal R}, \lambda_{\cal S}} < 1$
    and $\eta < 1$ is such that
    \begin{equation}\label{allodola}
        \eta^3
        \geq {
            \phi^*_{{\cal R}, \lambda_{\cal S}}
        \over
            \lambda
        }
        \qquad\hbox{and}\qquad
        \eta^4
        \geq \sqrt{
            \epsilon^*_{{\cal R}, \lambda_{\cal S}}
        }\,,
    \end{equation}
    then, setting
    $$
        \theta = {1 \over \eta^2}\left(
            {1 \over \lambda}
            \vee
            {
                \sqrt{\epsilon^*_{{\cal R}, \lambda_{\cal S}}}
            \over
                \eta \phi^*_{{\cal R}, \lambda_{\cal S}}
            }
        \right)
        \leq {\eta \over \phi^*_{{\cal R}, \lambda_{\cal S}}},
    $$
    it holds
    \begin{equation}\label{ben}
        \mathds{P}_{\mu_{\cal R}}\left(
            \theta < T_1,
            \sup_{t < T_1 - \theta} \left|
                A_\theta(t, f) - \mu_{\cal R}(f)
            \right| \leq 4 \eta \|f\|_\infty
        \right)
        \geq \left\{
            1 - 4 \eta - \sqrt{
                \epsilon^*_{{\cal R}, \lambda_{\cal S}}
            \over
                1 - \epsilon^*_{{\cal R}, \lambda_{\cal S}}
            }
        \right\}
    \end{equation}
    for all $f : {\cal X} \rightarrow \mathds{R}$.
\end{prp}

\par\noindent
{\it Proof:\/}
We can assume,
without loss of generality,
that $|f|$ is bounded by 1
and we will consider four events,
each of them with a small probability,
such that when none of them occurs
$T_1 > \theta$
and our time averages are indeed close
to $\mu_{\cal R}(f)$.
The first of these events
is $A = \{T^{\cal R}_{\lambda_{\cal S}} \leq \theta\}$.
When its complementary occurs
$\{\theta < T_1\}$ is implied.
By Proposition~\ref{monza}, it holds
\begin{equation}\label{scena}
    \mathds{P}_{\mu_{\cal R}}(A)
    \leq \mathds{P}_{\mu_{\cal R}}\left(
        T^{\cal R}_{\lambda_{\cal S}} \leq \theta
    \right)
    \leq {1 \over 2} \sqrt{
        \epsilon^*_{{\cal R}, \lambda_{\cal S}}
        \over 1 - \epsilon^*_{{\cal R}, \lambda_{\cal S}}
    } + 1 - e^{-\phi^*_{{\cal R}, \lambda_{\cal S}} \theta}
    \leq {1 \over 2} \sqrt{
        \epsilon^*_{{\cal R}, \lambda_{\cal S}}
        \over 1 - \epsilon^*_{{\cal R}, \lambda_{\cal S}}
    } + \phi^*_{{\cal R}, \lambda_{\cal S}} \theta.
\end{equation}

Our next event is $B = \{\ell_{\cal S}(T_1) > \theta'\}$
with
$$
    \theta' = \eta \theta
$$
being a small fraction of $\theta$.
When its complementary occurs
the excursions of $X$ in ${\cal R} \backslash {\cal S}$
up to $T_1$,
which cannot be longer than $\ell_{\cal S}(T_1)$,
will have a negligible contribution to time averages
on time scale $\theta$.
Since $\ell_{\cal S}(T_1)$ is an exponential variable
of rate $\lambda$, it holds
\begin{equation}\label{palco}
    \mathds{P}_{\mu_{\cal R}}(B)
    = \mathds{P}_{\mu_{\cal R}}\left(
        \ell_{\cal S}(T_1) > \theta'
    \right) = e^{-\lambda \theta'}.
\end{equation}

To control the supremum appearing in \eqref{ben},
we will divide the local time interval
$[0, T^{\cal R}_{\lambda_{\cal S}}]$
into intervals of length $\theta'$,
or smaller as far as the last one is concerned.
For $k_1$ a parameter that we will fix later,
our third event is
$C = \{T^{\cal R}_{\lambda_{\cal S}} > k_1 \theta'\}$.
When its complementary occurs it is sufficient
to control the time averages on time scale $\theta'$
for at most $k_1$ intervals of this length
to estimate time averages on time scale $\theta$.
By Proposition~\ref{monza}, it holds
\begin{equation}\label{penna}
    \mathds{P}_{\mu_{\cal R}}(C)
    = \mathds{P}_{\mu_{\cal R}}\left(
        T^{\cal R}_{\lambda_{\cal S}} > k_1 \theta'
    \right)
    \leq {1 \over 2} \sqrt{
        \epsilon^*_{{\cal R}, \lambda_{\cal S}}
        \over 1 - \epsilon^*_{{\cal R}, \lambda_{\cal S}}
    } + e^{-\phi^*_{{\cal R}, \lambda_{\cal S}} k_1 \theta'}.
\end{equation}

Our last event, which we will call $D$,
is that there is a local time interval
$[k \theta', (k + 1) \theta']$ with $k < k_1$
for which the time average associated with
the trace $X^{\cal R} = X_{{\cal R}, 0_{\cal S}}$
of $X$ on ${\cal R}$,
$$
    A^{\cal R}_{\theta'}(k \theta')
    = {1 \over \theta'}
    \int_{k \theta'}^{(k + 1) \theta'}
    f(X^{\cal R}(s))\,ds,
$$
differs from $\mu_{\cal R}(f)$
of more the $\eta$.
Any time interval $[t, t + \theta]$
with $t + \theta < T_1$
can be associated
``by removal of the excursions outside ${\cal R}$''
with a local time interval $[\ell_{\cal R}(t), \ell_{\cal R}(t + \theta)]$,
which, in turn, can be divided into at most
$\theta / \theta'$ fully covered
and two partially covered local time intervals of the form
$[k \theta', (k + 1) \theta']$.
Hence, if none of our last three events $B$--$D$ occurs,
then, for each $t < T_1 - \theta$
and using $|f| \leq 1$, it holds
$$
    \Bigl|
        A_\theta(t, f) - \mu_{\cal R}(f)
    \Bigr|
    \leq {\theta' + {\theta \over \theta'} \theta' \eta +  2 \theta' \over \theta}
     = 4 \eta.
$$
As far as the probability of $D$
is concerned, since $X^{\cal R}$ starts from equilibrium, the expected value
of $A^{\cal R}_{\theta'}(k\theta')$ is equal to $\mu_{\cal R}(f)$
and we can make an exact computation of its variance
by writing the decomposition $f$ in the eigenfunctions basis
associated with the generator of $X^{\cal R} = X_{{\cal R}, 0_{\cal S}}$.
By Proposition~\ref{pavia} its spectral gap,
which is associated with a zero value of $\lambda$,
is larger than $\gamma_{{\cal R}, \lambda_{\cal S}}$
and we get the upper bound
$$
    \mathds{E}_{\mu_{\cal R}}\left[
        \Bigl(
            A^{\cal R}_{\theta'}(k \theta')
            - \mu_{\cal R}(f)
        \Bigr)^2
    \right]
    \leq {
        2 {\rm Var}_{\mu_{\cal R}}(f)
        \over \gamma_{{\cal R}, \lambda_{\cal S}} \theta'
    }
    \leq {2 \over \gamma_{{\cal R}, \lambda_{\cal S}} \theta'}
    \,.
$$
We then get
\begin{equation}\label{maschera}
    \mathds{P}_{\mu_{\cal R}}(D)
    \leq {2 k_1 \over \gamma_{{\cal R}, \lambda_{\cal S}} \theta' \eta^2}
    \,.
\end{equation}

For \eqref{penna} and \eqref{maschera}
to be useful we need to choose $k_1$ in such a way
that, for $\eta \ll 1$,
$$
    k_1 \gg {1 \over \phi^*_{{\cal R}, \lambda_{\cal S}} \theta'}
    \qquad \hbox{and} \qquad
    k_1 \ll \gamma_{{\cal R}, \lambda_{\cal S}} \theta' \eta^2,
$$
which is possible if
\begin{equation}\label{charles}
    \theta'^2 \gg {
        1
        \over \eta^2 \gamma_{{\cal R}, \lambda_{\cal S}} \phi^*_{{\cal R}, \lambda_{\cal S}}
    }\,.
\end{equation}
We will choose $\theta'$ and $\eta$
to ensure this condition
and we will set
\begin{equation}\label{fave}
    k_1 = \sqrt{
        {1 \over \phi^*_{{\cal R}, \lambda_{\cal S}} \theta'}
        \gamma_{{\cal R}, \lambda_{\cal S}} \theta' \eta^2
    } = {\eta \over \sqrt{\epsilon^*_{{\cal R}, \lambda_{\cal S}}}}
    \,.
\end{equation}

For \eqref{scena} and \eqref{palco}
to be useful we need
to choose $\theta'$ in such a way that,
for $\eta \ll 1$,
(recall that $\theta = \theta' / \eta$)
$$
    \theta' \gg {1 \over \lambda}
    \qquad \hbox{and} \qquad
    \theta' \ll {\eta \over \phi^*_{{\cal R}, \lambda_{\cal S}}}
    \,.
$$
Once combined with \eqref{charles},
this reads
\begin{equation}\label{emilie}
    {1 \over \lambda}
    \vee {
        \sqrt{\epsilon^*_{{\cal R}, \lambda_{\cal S}}}
        \over \eta \phi^*_{{\cal R}, \lambda_{\cal S}}
    }
    \ll \theta'
    \ll {\eta \over \phi^*_{{\cal R}, \lambda_{\cal S}}}
    \,.
\end{equation}
We choose
$$
    \theta'
    = {1 \over \eta} \left(
        {1 \over \lambda}
        \vee
        {
            \sqrt{\epsilon^*_{{\cal R}, \lambda_{\cal S}}}
        \over
            \eta \phi^*_{{\cal R}, \lambda_{\cal S}}
        }
   \right)
$$
to satisfy the first inequality in~\eqref{emilie};
Condition~\eqref{allodola}
will ensure the last one.

Estimates \eqref{scena}--\eqref{maschera}
give then, recalling \eqref{fave} and using \eqref{allodola},
\begin{align*}
    &\mathds{P}_{\mu_{\cal R}}\left(
        T_1 \leq \theta \hbox{ or }
        \exists t < T_1 - \theta,\,
        \Bigl|A_\theta(t, f) - \mu_{\cal R}(f)\Bigr| > 4\eta
    \right)\\
    &\qquad
    \leq {1 \over 2} \sqrt{
        \epsilon^*_{{\cal R}, \lambda_{\cal S}}
    \over
        1 - \epsilon^*_{{\cal R}, \lambda_{\cal S}}
    } + {\phi^*_{{\cal R}, \lambda_{\cal S}} \over \eta^2}\left(
        {1 \over \lambda}
        \vee {
            \sqrt{\epsilon^*_{{\cal R}, \lambda_{\cal S}}}
        \over
            \eta \phi^*_{{\cal R}, \lambda_{\cal S}}
        }
    \right)
    + \exp\left\{
        - {\lambda \over \eta} \left(
            {1 \over \lambda}
            \vee {
                \sqrt{\epsilon^*_{{\cal R}, \lambda_{\cal S}}}
            \over
                \eta \phi^*_{{\cal R}, \lambda_{\cal S}}
            }
        \right)
    \right\}\\
    &\qquad\quad + {1 \over 2} \sqrt{
        \epsilon^*_{{\cal R}, \lambda_{\cal S}}
    \over
        1 - \epsilon^*_{{\cal R}, \lambda_{\cal S}}
    } + \exp\left\{
        - \phi^*_{{\cal R}, \lambda_{\cal S}} {
            \eta / \sqrt{\epsilon^*_{{\cal R}, \lambda_{\cal S}}}
        \over
            \eta
        }\left(
            {1 \over \lambda}
            \vee {
                \sqrt{\epsilon^*_{{\cal R}, \lambda_{\cal S}}}
            \over
                \eta \phi^*_{{\cal R}, \lambda_{\cal S}}
            }
        \right)
    \right\}
    + {
        2 \eta / \sqrt{\epsilon^*_{{\cal R}, \lambda_{\cal S}}}
    \over
        \gamma_{{\cal R}, \lambda_{\cal S}}
        {1 \over \eta} \left(
            {1 \over \lambda}
            \vee
            {
                \sqrt{\epsilon^*_{{\cal R}, \lambda_{\cal S}}}
            \over
                \eta \phi^*_{{\cal R}, \lambda_{\cal S}}
            }
        \right)
        \eta^2
    }\\
    &\qquad \leq \sqrt{
        \epsilon^*_{{\cal R}, \lambda_{\cal S}}
    \over
        1 - \epsilon^*_{{\cal R}, \lambda_{\cal S}}
    } + \eta + 2 e^{-1 / \eta} + 2 \eta
    \leq \sqrt{
        \epsilon^*_{{\cal R}, \lambda_{\cal S}}
    \over
        1 - \epsilon^*_{{\cal R}, \lambda_{\cal S}}
    } + 4\eta.
\end{align*}
\qed

\bibliographystyle{alpha}
\bibliography{bibliografia}

\end{document}